\documentclass[12pt]{article}

\usepackage{amssymb,amsmath,exscale}

\usepackage{graphicx}
\usepackage{dsfont}
\usepackage[applemac]{inputenc}
\usepackage[T1]{fontenc}

\usepackage{color}
\definecolor{marin}{rgb}   {0.,   0.3,   0.7} 
\definecolor{rouge}{rgb}   {0.8,   0.,   0.} 
\definecolor{sepia}{rgb}   {0.8,   0.5,   0.}

\usepackage[colorlinks,citecolor=marin,linkcolor=rouge,
            bookmarksopen,
            bookmarksnumbered
           ]{hyperref}

\newtheorem{lemma}{Lemma}[section]

\newtheorem{theorem}[lemma]{Theorem}
\newtheorem{proposition}[lemma]{Proposition}
\newtheorem{corollary}[lemma]{Corollary}
\newtheorem{remark}[lemma]{Remark}
\newtheorem{commentary}[lemma]{Commentary}
\newtheorem{example}[lemma]{Example}

\newtheorem{notation}[lemma]{Notation}
\newtheorem{definition}[lemma]{Definition}
\newtheorem{conclusion}[lemma]{Conclusion}
\numberwithin{equation}{section}
\newcommand{\QED}{\mbox{}\hfill \raisebox{-0.2pt}{\rule{5.6pt}{6pt}\rule{0pt}{0pt}} 
          \medskip\par}             
\newenvironment{Proof}{\noindent
    \parindent=0pt\abovedisplayskip = 0.5\abovedisplayskip
    \belowdisplayskip=\abovedisplayskip{\bfseries Proof. }}{\QED}
\newenvironment{Proofof}[1]{\noindent
    \parindent=0pt\abovedisplayskip = 0.5\abovedisplayskip
    \belowdisplayskip=\abovedisplayskip{\bfseries Proof of #1. }}{\QED}

\newcommand{\ad}{\mathrm{ad}}

\newcommand{\ab}{\boldsymbol{a}}
\newcommand{\bb}{\boldsymbol{b}}

\newcommand{\dd}{\mathrm{d}}
\newcommand{\deltab}{\boldsymbol{\delta}}

\newcommand{\Hc}{\mathcal{H}}

\newcommand{\M}{\mathcal{M}}
\newcommand{\Ic}{\mathcal{I}}

\newcommand{\jb}{{\boldsymbol{j}}}
\newcommand{\kb}{{\boldsymbol{k}}}
\newcommand{\Lc}{\mathcal{L}}

\newcommand{\ellb}{{\boldsymbol{\ell}}}
\newcommand{\N}{\mathbb{N}}
\newcommand{\Nc}{\mathcal{N}}

\newcommand{\Mc}{\mathcal{M}}
\newcommand{\Pc}{\mathcal{P}}
\newcommand{\Pcc}{\mathcal{P}^c}

\newcommand{\R}{\mathbb{R}}

\newcommand{\C}{\mathbb{C}}

\newcommand{\T}{\mathbb{T}}
\newcommand{\Z}{\mathbb{Z}}

\newcommand{\Sc}{\mathcal{S}}

\newcommand{\Zc}{\mathcal{Z}}
\newcommand{\xib}{\boldsymbol{\xi}}

\newcommand{\taub}{\boldsymbol{\tau}}

\newcommand{\Norm}[2]{\|#1\|\left.\vphantom{T_{j_0}^0}\!\!\right._{#2}}         
             
\title{Hamiltonian interpolation of splitting approximations for nonlinear  PDEs }        
\author{Erwan Faou{\small$\,^1$} and Beno\^it Gr\'ebert{\small$\,^2$}\\[4ex]
{\small$\,^1$}\, \small INRIA \& ENS Cachan Bretagne,  \\[-1ex]
\small Avenue Robert Schumann F-35170 Bruz, France. \\[-1ex]
\it \small email: \tt Erwan.Faou@inria.fr\\[2ex]
{\small$\,^2$}\, \small Laboratoire de Math\'ematiques Jean Leray,
Universit\'e de Nantes,\\[-1ex]
\small 2, rue de la Houssini\`ere
F-44322 Nantes cedex 3, France. \\[-1ex]
\small\it email: \tt benoit.grebert@univ-nantes.fr
}  

\begin{document}

\maketitle

\abstract{
We consider a wide class of semi linear Hamiltonian partial differential equations and their approximation by time splitting methods. We assume that the nonlinearity is polynomial, and that the numerical trajectory remains  at least  uniformly integrable with respect to an eigenbasis of the linear operator (typically the Fourier basis). We show the existence of a modified interpolated Hamiltonian equation whose exact solution coincides with the discrete flow at each time step over a long time depending on a non resonance condition satisfied by the stepsize. We introduce a class of modified splitting schemes fulfilling this condition at a high order and prove for them that the numerical flow and the continuous flow remain close over exponentially long time with respect to the step size. For standard splitting or implicit-explicit scheme, such a backward error analysis result holds true on a time depending on a cut-off condition in the high frequencies (CFL condition). This analysis is valid in the case where the linear operator has a discrete (bounded domain) or continuous (the whole space) spectrum. \\[2ex]
}\\[2ex]


\section{Introduction}

The {\em Hamiltonian interpolation} of a symplectic map which is a perturbation of the identity is a central problem both in the study of Hamiltonian systems and in their discretization by numerical methods.  This question actually goes back to Moser \cite{Moser68} who interpreted such a map as the exact flow of a Hamiltonian system in the  finite dimensional context. Such a result was later refined and extended by {Benettin and Giorgilli} \cite{BeGi94} to the analysis of symplectic numerical methods, and leads to the seminal backward error analysis results of {Hairer, Lubich} \cite{HL97} and {Reich} \cite{Reich99} for numerical integrators applied to ordinary differential equations. These results constitute now a cornerstone of the modern geometric numerical integration theory (see \cite{HLW,Reic04}). 

In the finite dimensional case, the situation can be described as follows: if $(p,q)\mapsto \Psi^h(p,q)$ is a symplectic map from  the phase space $\R^{2d}$ into itself, and if this map is a perturbation of the identity, $\Psi^h \simeq {\rm Id} + \mathcal{O}(h)$, then there exists a Hamiltonian function $H_h(p,q)$ such that $\Psi^h$ can be interpreted as the flow at time $t = h$ of the Hamiltonian system associated with $H_h$. If $\Psi^h$ is analytic, such results holds up to an error that is exponentially small with respect to the small parameter $h$, and for $(p,q)$ in a compact set of the phase space. In application to numerical analysis, the map $\Psi^h(p,q) \simeq \Phi^h_H(p,q)$ is a numerical approximation of the exact flow associated with a given Hamiltonian $H$. As a consequence, the modified Hamiltonian $H_h$, which turns out to be a perturbation of the initial Hamiltonian $H$, is preserved along the numerical solution, $(p_n,q_n) = (\Psi^h)^{n}(p_0,q_0)$, $n \in \N$. More precisely we have 
\begin{equation}
\label{Efinited}
(p_n,q_n) = (\Psi^h)^{n}(p_0,q_0) = (\Phi^{nh}_{H_h}) (p_0,q_0) + nh \exp(-1/(ch))
\end{equation}
under the a priori assumption that the sequence $(p_n,q_n)_{n \in \N}$ remains in a compact set used to derive the analytic estimates. Here the constant $c$ depends on the eigenvalues of the quadratic part of the Hamiltonian $H$ (i.e.\ the linear part of the associated ODE). As a consequence the qualitative behavior of the discrete dynamics associated with the map $\Psi^h$ over exponentially long time can be pretty well understood through the analysis of the continuous system associated with $H_h$.  

The extension of such results to Hamiltonian partial differential equations (PDE) faces the principal difficulty that the Hamiltonian function involves operators with unbounded eigenvalues making the constant $c$ in the previous estimate blow up. 
The goal of this work is to overcome this difficulty and to give Hamiltonian interpolation results for splitting methods applied to semi linear Hamiltonian PDEs. 

We consider a class of Hamiltonian PDEs associated with a Hamiltonian $H$ than can be split into 
a quadratic functional associated with an unbounded linear operator, $H_0$, and a nonlinearity $P$ which is a   polynomial  functional at least cubic:
$$
H = H_0 + P.
$$
Typical examples are given by the non linear Schr\"odinger equation (NLS)
\begin{equation}
\label{NLS1}
i \partial_t u = \Delta u + f(u,\bar u)
\end{equation}
or the non linear wave equation (NLW)
\begin{equation}
\label{NLW}
\partial_{tt} u - \Delta u = g(u)
\end{equation}
both set on the torus $\T^d$ or on the whole space $\R^d$. Here $f$ and $g$ are   polynomials  having a zero of order at least three at the origin, for instance $f(u,\bar u)=|u|^2u$ for the cubic defocusing (NLS) and $g(u)=u^3$ for the classical (NLW).\\
 To approximate such equations, splitting methods are widely used: They consist in decomposing the exact flow $\Phi_{H}^h$ at a small time step $h$ as compositions of the flows $\Phi_{H_0}^h$ and $\Phi^h_P$. The Lie-Trotter splitting methods
\begin{equation}
\label{Eclasssplit}
u(nh,x)\sim u_n(x):=( \Phi_{H_0}^h \circ \Phi^h_P)^nu_0(x) \quad \mbox{or} \quad u_n(x):=(  \Phi^h_P\circ \Phi_{H_0}^h)^nu_0(x)
\end{equation}
are known to be order 1 approximation scheme in time, when the solution is {\em smooth}.
Here {\em smooth} means that the numerical solution belongs to some Sobolev space $H^s$, $s > 1$, uniformly in time, see \cite{Jahn00,Hans08} for the linear case and \cite{L08a} for the non linear Schr\"odinger equation.

These schemes are all symplectic, preserve the $L^2$ norm if $H_0$ and $P$ do, and can be easily implemented in practice. For instance, in the cases of (NLS) or (NLW) just above,  we can diagonalize the linear part  and integrate it by using fast Fourier transform and integrate the nonlinear part explicitly (it is an ordinary differential equation). More generally, instead of  fast Fourier transform, pseudo spectral methods are used if the spectrum of $H_0$ is known and available. 

\medskip
%

The Hamiltonian interpolation problem described above can be formulated here as follows: is it possible to find a {\em modified energy} (or modified Hamiltonian function) $H_h$ depending on $H_0$, $P$ and on the chosen stepsize, such that 
 \begin{equation}
\label{EBCHF}
\Phi^h_P\circ \Phi_{H_0}^h \simeq \Phi_{H_h}^h.
\end{equation}
Formally, this question corresponds to the classical Baker-Campbell-Hausdorff (BCH) formula (see \cite{Baker05,Haus06}) for which the modified Hamiltonian $H_h$ 
expresses as iterated Poisson brackets between the two Hamiltonian $P$ and $H_0$. Hence we observe that the validity of such representation {\em a priori} depends on the smoothness of the discrete solution $u_{n}$. But this is not fair, as there is no reason for $u_n$ to be, a priori, uniformly smooth over long time. Actually, it is well known that this assumption is not satisfied in general - even if the exact solution of the continuous system is smooth - making such analysis irrelevant in many cases of applications, see for instance \cite{FG2}. Note that this smoothness assumption balances the blow up of the constant $c$ in the finite dimensional exponential estimates \eqref{Efinited}. In other words, the norm for which \eqref{EBCHF} is valid with the ``finite dimensional'' method is not the norm required for the analytic estimates of the error in the infinite dimensional context. 

\medskip 

In this work, we give backward error analysis results in the spirit of \cite{BeGi94,HL97,Reich99} for the map \eqref{EBCHF}. We follow a general approach recently  developed in \cite{DF09} for linear PDEs.  
In this linear context, Debussche and Faou proved a formula of the form \eqref{EBCHF} by considering smoothed schemes, namely schemes where $\Phi_{H_0}^h$ is replaced by the midpoint approximation of the unbounded part $H_0$. This yields schemes of the form
\begin{equation}
\label{Edefsplit}
\Phi^h_P\circ \Phi_{A_0}^1 \quad \mbox{or} \quad  \Phi^1_{A_0}\circ \Phi_{P}^h
\end{equation}
where $A_0$ depends on $h$, and is an approximation of the operator $ h H_0$ when $h \to 0$. 

In the present work we consider the context of nonlinear PDEs ($P$ is at least cubic) and we assume that the eigenvalues of $A_0$ are of the form
 \begin{equation}
\label{Elamb}
\lambda_a = \alpha_h(h \omega_a)
\end{equation}
where  we denote by $\omega_a$ the frequencies of $H_0$ (here the index $a$ belongs to a discrete set or is a continuous variable), and
where $\alpha_h$ is a {\em filter} function satisfying $\alpha_h(x) \simeq x$ for small $x$. It turns out that under such an assumption, \eqref{Edefsplit} remains an order one approximation of the continuous solution, provided that the numerical solution is smooth. This general setting contains the {\em standard} splitting methods \eqref{Eclasssplit} which corresponds to the choice $\alpha_h(x) = x$ for all $x$, and  the midpoint approximation of the linear flow  which corresponds to the function  (see \cite{DF09})
$$
\alpha_h(x) = 2 \arctan (x/2). 
$$


\medskip

Our results
can be described as follows:  we  assume that, in the Fourier variables associated with the diagonalization in an eigenbasis of $H_0$, the numerical trajectory remains bounded  in the space $\ell^1_s$, for some fixed $ s \geq 0$.  In the case of the Laplace operator on $\T^d$, this means to assume that $u_n$ has its Fourier transform satisfying $(|k|^s \widehat{u_n}(k)_{k\in \Z^d})$  in $\ell^1(\Z^d)$ while  in the case of the same operator on $\R^d$, this corresponds to assume that $u_n$ has its Fourier transform satisfying  $|\xi|^s \widehat{u_n}(\xi)$  in $L^1_{\xi}(\R^d)$. When for instance $s = 0$, this implies in both cases that the numerical trajectory  remains  bounded in the $L^\infty$ space  associated with the space variable uniformly in time.  This is actually the standard assumption is the finite dimensional case: no blow up in finite time.\\
 Then for $z$ (the Fourier variable) in a fixed ball of $\ell^1_s$,  we  construct a Hamiltonian $H_h$ such that 
\begin{equation}
\label{Eest1}
\Norm{\Phi_P^h \circ \Phi_{A_0}^1(z) - \Phi_{H_h}^h(z) }{\ell^1_s} \leq h^{N+1} (CN )^N
\end{equation}
for some constant $C$ depending on $P$  and $s$. 
It turns out that the number $N$ depends on a sort of non resonance condition satisfied by the eigenvalues of the operator $A_0$. With $\lambda_a$ the eigenvalues of $A_0$ (see \eqref{Elamb}) this condition can be written
\begin{equation}
\label{Enonresres}
\forall\, j = 1,\ldots, r, \quad \forall\, (a_1,\ldots,a_j)\quad | \lambda_{a_1} \pm \cdots \pm \lambda_{a_j} | < 2\pi.
\end{equation}
If this is satisfied for some $r\geq 2$, then we can take  in the previous estimate
$$
N = (r-2)/(r_0 - 2)
$$ 
 where $r_0$ is the degree of the polynomial $P$.  Note that this condition is independent of the regularity index $s$. We thus see that the actual existence of a modified PDEs almost coinciding with the numerical scheme at each time step depends on the degree of the polynomial $P$ and on the non resonance condition \eqref{Enonresres}.

For numerical schemes associated with a filter function 
\begin{equation}
\label{Egood}
\alpha_h(x) = \sqrt{h} \arctan( x / \sqrt{h}),
\end{equation}
 we can actually prove that the non resonance condition \eqref{Enonresres} is satisfied for large $r \simeq 1/\sqrt{h}$. The analytic estimate \eqref{Eest1} then yields an exponentially small error at each step. The choice \eqref{Egood} typically induces a stronger regularization in the high frequencies than the midpoint rule, without breaking the order of approximation of the method. 

For the other classical schemes (standard splitting, implicit-explicit schemes), the non resonance condition \eqref{Enonresres} is in general not satisfied unless a Courant-Friedrichs-Lewy (CFL) condition is imposed \cite{CFL28}, depending on the wished number $r$. Although   the  CFL number tends to $0$ when $r $ goes to $ \infty$ for all schemes, 
 we obtain very different results 
for relatively small $r$ regarding on the method and on the specific form of the nonlinearity. For example the combination of the midpoint approximation of the linear part in an implicit-explicit scheme together with a CFL condition with CFL number of order $1$ in the case of the cubic nonlinear Schr\"odinger equation yields a modified equation with precision of order $h^7$ at each step (see section 5 for details and other examples). 
\medskip


We can precise in which sense the estimate \eqref{Eest1} induces a better control of the numerical solution. Actually,  the modified Hamiltonian depends on $A_0$ and $P$ and can roughly speaking be described as a Hamiltonian $A_0/h + \widetilde{P}$ where $\widetilde{P}$ is a modified nonlinear Hamiltonian. For example in the case of the method  \eqref{Egood} applied to the cubic  (NLS)  on $\R^d$  we  show (cf. section 6) that  the modified energy reads
\begin{equation}
\label{EmodEN}
\int_{\R^d}  \frac{1}{\sqrt{h}}\arctan(\sqrt{h} |\xi|^2 ) |\widehat{u}(\xi)|^2 \dd \xi + P_h^{(1)}(u), 
\end{equation}
and is preserved up to an error $\mathcal{O}(h)$ over exponentially long times $t \leq e^{c/\sqrt{h}}$,  provided the numerical solution remains bounded in $\ell^1 = \ell^1_0$ only.  Here, $\widehat{u}$ denotes the Fourier transform of $u$ and  $P_h^{(1)}(u)$ is a polynomial nonlinearity of degree $4$.  The actual expression of  $P_h^{(1)}(u)$   depends on small denominators controlled by the non resonance condition \eqref{Enonres} and is given in Theorem \ref{t:nls}. Note that this modified energy will be close to the original energy only if $u$ is smooth, something that is not guaranteed by the present analysis  unless the a priori regularity index $s$ is large enough. However, the preservation of \eqref{EmodEN} implies the control of the $H^1$ norm of low modes of the numerical solution, and  of  the $L^2$ norm of high modes, as in 
\cite{DF09}, over very long time. 


\medskip

The relative question of persistence of smoothness of the numerical solution has recently known many progresses: see  \cite{DF07,FGP1,FGP2,CHL08b,CHL08c,GL08a,GL08b} and \cite{FG2}. However, we emphasize that the  resonances considered in all these works are not of the same type as here. They read
\begin{equation}\label{nonresk}
 \lambda_{a_1} \pm \cdots \pm \lambda_{a_r} \simeq 2k\pi \quad \mbox{ for some } k\in \Z
 \end{equation}
and  they control the persistence of the regularity of the solution over long time. As explain in  \cite{FGP1,FGP2, FG2}, the case $k=0$ corresponds to {\em physical resonances} (they actually are  close to  resonances of the PDE under study, see also  \cite{BG06,Greb07}) while the case  $k\neq 0$ corresponds to {\em numerical resonances} (they are not resonances of the PDE).\\
  In the present work, we address a different question and we point out that  possible physical resonances do not contradict the existence of the modified energy \eqref{EmodEN}. Such relations are in fact also resonances of the  modified equation associated with the modified Hamiltonian $H_h$ . In such a case, no preservation of the regularity in high index  Sobolev spaces $H^m$  for the numerical solution can be expected, though the existence of the modified equation is guaranteed by the present analysis  even if $s = 0$. Only numerical resonances can avoid the existence of a modified energy.
  
 \medskip 
 
 We end this introduction by two important  commentaries:
\begin{commentary}
We emphasize that we use the assumption that the nonlinearity has {\em zero momentum}. This ensures the fact that the Poisson bracket of two Hamiltonian  with bounded polynomial coefficients  will still act on  the space  $\ell_s^1$, and that the corresponding vector field will be well defined on $\ell_s^1$.  This kind of assumption is made in \cite{FGP2} and precludes Hamiltonian PDEs with strong non local nonlinearity. Note however that we may obtain results with small momentum by combining the results in \cite{DF09}, i.e.\ by assuming that the nonlinearity $P(x,u)$ depends smoothly on $x$ and polynomially on $u$. 
\end{commentary}

\begin{commentary}
Note that similar results can be obtained for Strang-like splitting methods generalizing \eqref{Edefsplit}, i.e. schemes of the form 
$$
\Phi_P^{h/2} \circ \Phi_{A_0}^1 \circ \Phi_P^{h/2}
$$
for which the analysis is very close to the one developed here. Similarly, we conjecture that the same kind of results holds true for any symplectic Runge-Kutta method under CFL condition, but this leads to technical difficulties. The big advantage of the splitting approach is that the induction equations building the modified Hamiltonian directly  express at the level of the Hamiltonian, and not on the corresponding vector field (in particular we do not have to show that the vector field constructed at each step is - indeed - Hamiltonian). 
\end{commentary}

\section{Abstract Hamiltonian formalism}

In this section we describe the general form of Hamiltonian PDEs we are considering. We focus  on the case where the spectrum of the linear operator is discrete: elliptic operator on the torus or on a bounded domain with boundary conditions. We give in the last section of this paper the principal changes that have to be made to obtain similar results in the case where the spectrum of the linear operator is continuous (Schr\"odinger or wave equation over the whole space). 

So typically, the solution of the PDE under study is decomposed in the eigenbasis  of the linear part:

$$\psi(t,x)=\sum \xi_k(t) \phi_k(x)$$
and we observe the PDE in the Fourier-like variables $\xi=(\xi_k)$.

The main difference with the presentation in \cite{FGP1,FGP2,Greb07} lies in the choice of the phase space. We consider Fourier variables $\xi$ belonging to $\ell_s^1$ and not to $\ell^2_s$ the weighted $\ell^2$ space which corresponds to function $\psi$ belonging to the Sobolev spaces $H^s$.  In particular,  considering  Fourier coefficients in  $\ell^1_0$  remains to consider bounded functions, i.e. to assume that the solution (or its numerical approximation) of the PDE under study is (essentially)  bounded.   

\subsection{Setting and notations}

We set  $\Nc = \Z^d$ or $\N^d$ (depending on the concrete application) for some $d \geq 1$.  For $a = (a_1,\ldots,a_d) \in \Nc$, we set
$$
|a|^2 = \max\big(1,a_1^2 + \cdots + a_d^2\big). 
$$
We consider the set of variables $(\xi_a,\eta_b) \in \C^{\Nc} \times \C^{\Nc}$ equipped with the symplectic structure
\begin{equation}
\label{Esymp}
i \sum_{a \in \Nc} \dd \xi_a \wedge \dd \eta_a. 
\end{equation}
We define the set $\Zc = \Nc \times \{ \pm 1\}$. For $j = (a,\delta) \in \Zc$, we define $|j| = |a|$ and we denote by $\overline{j}$ the index $(a,-\delta)$. 

We will identify a couple $(\xi,\eta)\in \C^{\Nc} \times \C^{\Nc}$ with 
$(z_j)_{j \in \Zc} \in \C^{\Zc}$ via the formula
$$
j = (a,\delta) \in \Zc  \Longrightarrow 
\left\{
\begin{array}{rcll}
z_{j} &=& \xi_{a}& \mbox{if}\quad \delta = 1,\\[1ex]
z_j &=& \eta_a & \mbox{if}\quad \delta = - 1.
\end{array}
\right.
$$
By  a slight abuse of notation, we often write $z = (\xi,\eta)$ to denote such an element. 
\begin{example}\label{ex} In the case where $H_0=-\Delta$ on the torus $\T^d$, the eigenbasis is the Fourier basis, $\Nc=\Z^d$ and $\xi$ is the sequence associated with a function $\psi$ while $\eta$ is the Fourier sequence associated with a function $\phi$ via the formula
$$\psi(x)=\sum_{a\in \Nc}\xi_a e^{ia.x}\quad \mbox{and} \quad \phi(x)=\sum_{a\in \Nc}\eta_a e^{-ia.x}.$$
\end{example}
 For  a given $s \geq 0$, 
we consider the Banach space $\ell_s^1 := \ell_s^1(\Zc,\C)$ made of elements $z \in \C^{\Zc}$ such that
$$
\Norm{z}{\ell_s^1} := \sum_{j \in \Zc} |j|^s |z_j| < \infty,
$$
and equipped with the symplectic form \eqref{Esymp}. We will often write simply $\ell^1 = \ell^1_0$.  We moreover define for $s > 1$ the Sobolev norms 
$$
\Norm{z}{H^s} = \Big(\sum_{j \in \Zc} |j|^{2s} |z_j|^2\Big)^{1/2}. 
$$

For a function $F$ of $\mathcal{C}^1(\ell^1_s,\C)$, we define its gradient by 
$$
\nabla F(z) = \left( \frac{\partial F}{\partial z_j}\right)_{j \in \Zc}
$$
where by definition, we set for $j = (a,\delta) \in \Nc \times \{ \pm 1\}$, 
$$
 \frac{\partial F}{\partial z_j} =
  \left\{\begin{array}{rll}
 \displaystyle  \frac{\partial F}{\partial \xi_a} & \mbox{if}\quad\delta = 1,\\[2ex]
 \displaystyle \frac{\partial F}{\partial \eta_a} & \mbox{if}\quad\delta = - 1.
 \end{array}
 \right.
$$
Let $H(z)$ be a function defined on $\ell^1_s$. If $H$ is smooth enough, we can associate with this function the Hamiltonian vector field $X_H(z)$ defined by
$$
X_H(z) = J \nabla H(z) 
$$
where $J$ is the symplectic operator on $\ell^1_s$ induced  by the symplectic form \eqref{Esymp}.

For two functions $F$ and $G$, the Poisson Bracket is (formally) defined as
\begin{equation}\label{poisson}
\{F,G\} = \nabla F^T J \nabla G = i \sum_{a \in \Nc} \frac{\partial F}{\partial \eta_j}\frac{\partial G}{\partial \xi_j} -  \frac{\partial F}{\partial \xi_j}\frac{\partial G}{\partial \eta_j}.  
\end{equation}

We say that $z\in \ell^1_s$ is {\em real} when $z_{\overline{j}} = \overline{z_j}$ for any $j\in \Zc$. In this case, we write $z=(\xi,\bar\xi)$ for some $\xi\in \C^{\Nc}$. Further we say that a Hamiltonian function $H$ is 
 {\em real } if $H(z)$ is real for all real $z$. 
 
\begin{example} \label{ex2}   Following Example \ref{ex}, a real $z$ corresponds to the relation $\phi=\bar \psi$ and a typical real Hamiltonian reads $H(z)= \int_{\T^d}h(\psi(x),\phi(x))\dd x$ where $h$ is a regular function from $\C^2$ to $\C$ satisfying $h(\zeta,\bar \zeta)\in \R$ for all $\zeta\in \C$.
\end{example}

\begin{definition}\label{def:2.1}
 For a given $s \geq 0$, we  denote by $\Hc_s$ the space of real Hamiltonians $P$ satisfying 
$$
P \in \mathcal{C}^{1}(\ell^1_s,\C), \quad \mbox{and}\quad 
X_P \in \mathcal{C}^{1}(\ell^1_s,\ell^1_s). 
$$
\end{definition}
Notice that for $F$ and $G$ in $\Hc_s$ the formula \eqref{poisson} is well defined. 
\begin{remark}\label{remH0}
At this stage we have to note that, in general, the quadratic Hamiltonian $H_0$ corresponding to the linear part of our PDE will not be defined for $\xi\in \ell^1_s$. For instance for the NLS equation, $H_0$ is associated with the Laplace operator (see example \ref{ex} above) and reads in Fourier variable $\sum k^2 |\xi_k|^2$. Nevertheless it generates a flow which maps $\ell^1_s$ into $\ell^1_s$ and which is given for all time $t$ and for all indices $k$ by  $\xi_k(t)=e^{-ikt}\xi_k(0)$.
\end{remark}

We will verify later (see example \ref{ex3}) that  the typical real Hamiltonians given in Example \ref{ex2}  belong to the class $\Hc_s$. Actually the proof is not totally trivial  because the Fourier transform is not well adapted to the $\ell^1_s$ space. 

With a given Hamiltonian function $H \in \Hc_s$, we associate the Hamiltonian system
$$
\dot z = J \nabla H(z)
$$
which can be written
\begin{equation}
\label{Eham2}
\left\{
\begin{array}{rcll}
\dot\xi_a &=& \displaystyle - i \frac{\partial H}{\partial \eta_a}(\xi,\eta) & a \in \Nc,\\[2ex]
\dot\eta_a &=& \displaystyle i \frac{\partial H}{\partial \xi_a}(\xi,\eta)& a \in \Nc.  
\end{array}
\right.
\end{equation}
In this situation, we define the flow $\Phi_H^t(z)$ associated with the previous system (for an interval of times $t \geq 0$ depending a priori on the initial condition $z$). Note that if $z = (\xi,\bar \xi)$ and if  $H$ is real, the flow $(\xi^t,\eta^t) = \Phi_H^t(z)$ is also real for all time $t$ where the flow is defined: $\xi^t = \bar {\eta}^t$.
When $H$ is real, it may be useful to introduce the real variables $p_a$ and $q_a$ given by
$$
\xi_a = \frac{1}{\sqrt{2}} (p_a + i q_a)\quad \mbox{and}\quad \bar{\xi}_a =  \frac{1}{\sqrt{2}} (p_a - i q_a),
$$
the system \eqref{Eham2} is then equivalent to the system
$$
\left\{
\begin{array}{rcll}
\dot p_a &=& \displaystyle -  \frac{\partial H}{\partial q_a}(q,p) & a \in \Nc,\\[2ex]
\dot q_a &=& \displaystyle  \frac{\partial H}{\partial p_a}(q,p),& 	a \in \Nc,  
\end{array}
\right.
$$
where by a slight abuse of notation we still denote the Hamiltonian with the same letter: $H(q,p) = H(\xi,\bar\xi)$. 

%

We now describe the hypothesis needed on the Hamiltonian nonlinearity $P$.

Let $\ell \geq 2$. We consider $ \jb = (j_1,\ldots,j_\ell) \in \Zc^{\ell}$, and we set for all $i = 1,\ldots \ell$
$j_i = (a_i,\delta_i)$ where $a_i \in \Nc$ and $\delta_i \in \{\pm 1\}$. We define
$$
\overline \jb = (\overline{j}_1,\ldots,\overline j_\ell)\quad\mbox{with}\quad \overline{j}_i = (a_i,-\delta_i), \quad i = 1,\ldots,\ell.
$$
We also 
use the notation 
$$
z_\jb = z_{j_1}\cdots z_{j_\ell}. 
$$
We define the momentum $\Mc (\jb)$ of  the multi-index $\jb$ by
\begin{equation}
\label{EMcb}
\Mc(\jb) = a_{1} \delta_{1} + \cdots + a_\ell \delta_\ell.
\end{equation}
We then define the set of indices with zero momentum
\begin{equation}
\label{EIr}
 \Ic_\ell =  \{  \jb = (j_1,\ldots,j_\ell) \in \Zc^{\ell}, \quad \mbox{with}\quad \Mc(\jb) = 0\}.
\end{equation}
We can now define precisely the type of polynomial nonlinearities we consider:
\begin{definition}\label{zero}
We say that a polynomial Hamiltonian $P \in \Pc_k$ if $P$ is real, of degree $k$, have a zero of order at least $2$ in $z = 0$, and if 
\begin{itemize}
\item $P$ contains only monomials $a_\jb z_\jb$ having zero momentum, i.e.\ such that $\Mc(\jb)=0$ when $a_\jb \neq 0$ and thus $P$ formally reads
\begin{equation}
\label{EexpP}
P(z) = \sum_{\ell = 2}^k \sum_{\jb \in \Ic_\ell} a_{\jb} z_{\jb}
\end{equation}
with the relation $a_{\bar\jb} = \bar{a}_{\jb}$. 
\item The coefficients $a_{\jb}$ are bounded, i.e.\ satisfy 
$$
\forall\, \ell = 2,\ldots,k,\quad  \forall\, \jb = (j_1,\cdots,j_\ell) \in \Ic_\ell, \quad |a_\jb| \leq C. 
$$
In the following, we set
\begin{equation}
\label{EnormP}
\Norm{P}{} = \sum_{\ell = 2}^k \sup_{\jb \in \Ic_\ell} |a_\jb|. 
\end{equation} 
\end{itemize}
\end{definition}

\begin{definition}
We say that $P \in \Sc\Pc_k$ if $P \in \Pc_k$ has coefficients $a_{\jb}$ such that $a_{\jb} \neq 0$ implies that $\jb$ contains the same numbers of positive and negative indices: 
$$
\sharp\{i \,|\,  j_{i} = (a_i, +1)\} = \sharp\{i \,|\,  j_{i} = (a_i, - 1)\}.
$$
In other words, $P$ contains only monomials with the same numbers of $\xi_i$ and $\eta_i$. Note that this implies that $k$ is even. 
\end{definition}

\begin{example}\label{ex3}
Following example \ref{ex2}, $P(z)=\int_{\T^d}p(\psi(x),\phi(x))\dd x$, where $p$ is a polynomial of degree $k$ in $\C[X,Y]$ satisfying $p(\zeta,\bar \zeta)\in \R$ and having a zero of order at least 2 at the origin, defines  a Hamiltonian in $\Pc_k$. \\
An example of  polynomial Hamiltonian in $\Sc\Pc_{2k}$ is given by $P=\int |\psi|^{2k}\dd x$.
\end{example}

The zero momentum assumption  in Definition \ref{zero} is crucial in order to obtain the following Proposition:

\begin{proposition}
\label{P1}
Let $k \geq 2$  and $s \geq 0$, then we have $\Pc_k \subset \Hc_s$,  and for $P \in \Pc_k$,  we have the estimates
\begin{equation}
\label{Epot}
|P(z)| \leq \Norm{P}{}\big(\max_{n = 2,\ldots,k}  \Norm{z}{\ell^1_s}^n \big)
\end{equation}
and 
\begin{equation}
\label{Echamp}
\forall\, z \in \ell^1_s, \quad  \Norm{X_P(z) }{\ell^1_s} \leq 2 k (k-1)^{s} \Norm{P}{}\Norm{z}{\ell^1_s}\big(\max_{n = 1,\ldots,k-2}  \Norm{z}{\ell^1}^n \big). 
\end{equation}
Moreover, for $z$ and $y$ in $\ell^1_s$, we have 
\begin{equation}
\label{Elip}
\Norm{X_P(z) - X_P(y) }{\ell^1_s} \leq 4 k (k-1)^s \Norm{P}{}{} \big(\max_{n = 1,\ldots,k-2} (\Norm{y}{\ell^1_s}^n, \Norm{z}{\ell^1_s}^n )\big)\Norm{z-y}{\ell^1_s}.
\end{equation}
Eventually,  for $P\in \Pc_k$ and $Q \in \Pc_\ell$, then $\{P,Q\} \in \Pc_{k+ \ell - 2}$ and we have the estimate
\begin{equation}
\label{Ebrack}
\Norm{\{P,Q\}}{} \leq 2 k\ell\Norm{P}{}\Norm{Q}{}.
\end{equation}
If now $P \in \Sc\Pc_k$ and $Q \in \Sc\Pc_k$, then $\{P,Q\} \in \Sc\Pc_{k+ \ell - 2}$. 
\end{proposition}


\begin{remark} The estimate \eqref{Echamp} is a sort of tame estimate. Notice that the same estimate with $\ell^1_s$ replaced by $\ell^2_s$ is proved in \cite{Greb07}  under a decreasing assumption on the coefficients of the polynomial $P$. Actually  the present proof is much simpler. We also notice that, with this estimate and following \cite{Greb07}, we could develop a  Birkhoff normal form theory in $\ell^1_s$. The only disagreement with  this choice of Fourier space is that $\ell^1$ is not the image by the Fourier transform of $L^\infty$, making difficult the pull back of the normal form to the original variables.
\end{remark}

\begin{Proof}
Assume that $P$ is given by \eqref{EexpP}, and denote by $P_i$ the homogeneous component of degree $i$ of $P$, i.e. 
$$
P_i(z) = \sum_{\jb \in \Ic_i} a_\jb z_\jb, \quad i = 2,\ldots,k. 
$$
We have for all $z$
$$
|P_i(z)| \leq \Norm{P_i}{} \Norm{z}{\ell^1}^i  \leq \Norm{P_i}{} \Norm{z}{\ell_s^1}^i. 
$$
The first inequality \eqref{Epot} is then a consequence of the fact that 
\begin{equation}
\label{Epotpot}
\Norm{P}{} = \sum_{i = 2}^k \Norm{P_i}{}. 
\end{equation}
Now let $j = (a,\epsilon)\in \Zc$ be fixed. The derivative of a given monomial $z_{\jb} = z_{j_1}\cdots z_{j_i}$ with respect to $z_j$ vanishes except if $j \subset \jb$. Assume for instance that $j = j_i$. Then the zero momentum  condition implies that $\Mc(j_1,\ldots,j_{i-1}) = - \epsilon a$ and we can write 
\begin{equation}
\label{Emomentj}
|j|^s \left|\frac{\partial P_i}{\partial z_j} \right|  \leq i \Norm{P_i}{} \sum_{\jb \in \Zc^{i-1} \,  , \, \Mc(\jb) = - \epsilon a }  |j|^s |z_{j_1} \cdots z_{j_{i-1}}|.
\end{equation}
Now in this formula, for a fixed multiindex $\jb$, the zero momentum condition implies that 
\begin{equation}
\label{raslebol}
|j|^s \leq \big( |j_1| + \cdots + |j_{i-1}|\big)^s \leq (i-1)^s \max_{n = 1,\ldots,i-1} |j_n|^s. 
\end{equation}
Therefore, after summing in $a$ and $\epsilon$ we get
\begin{multline}
\label{takecare}
\Norm{X_{P_i}(z) }{\ell^1_s} \leq 2 i(i-1)^s \Norm{P_i}{} \sum_{\jb \in \Zc^{i-1}} \max_{n = 1,\ldots,i-1} |j_n|^s  |z_{j_1}| \cdots  | z_{j_{i-1}}| \\
\leq 2 i (i-1)^s \Norm{P_i}{}\Norm{z}{\ell^1}\Norm{z}{\ell^1_s}^{i-2}
\end{multline}
which yields \eqref{Echamp} after summing in $i = 2,\ldots,k$. 

Now for $z$ and $y$ in $\ell^1_s$, we have with the previous notations
$$
|j|^s \left|\frac{\partial P_i}{\partial z_j}(z) -\frac{\partial P_i}{\partial z_j}(y) \right|
\leq  \sum_{q \in \Zc}  |j|^s\left|\int_{0}^1 \frac{\partial P_i}{\partial z_j\partial z_q}(ty + (1-t)z) \dd t\right| |z_q - y_q|.
$$
But we have for fixed $j = (\epsilon,a)$ and $q = (\delta,b)$ in $\Zc$, and for all $u \in \ell^1_s$
$$
 |j|^s\left|\frac{\partial P_i}{\partial z_j\partial z_q}(u)\right| \leq  i \Norm{P_i}{} \sum_{\jb \in \Zc^{i-2} \,  , \, \Mc(\jb) = - \epsilon a- \delta b }  |j|^s |u_{j_1} \cdots u_{j_{i-2}}|.
$$
In the previous sum, we necessarily have that $\Mc(\jb,j,k) = 0$, and hence 
$$
|j|^s \leq \big( |j_1| + \cdots + |j_{i-2}| + |q|\big)^s 
\leq (i-1)^s |q|^s \prod_{n = 1}^{i-2}|j_n|^s. 
$$
Let $u(t) = ty + (1-t)z$, we have for all $t \in [0,1]$ with the previous estimates
\begin{multline*}
 |j|^s\left|\int_{0}^1\frac{\partial P_i}{\partial z_j\partial z_q}(u(t)) \dd t\right| \leq \\
i (i-1)^s |q|^s  \Norm{P_i}{}  \int_0^1 \sum_{\jb \in \Zc^{i-2} \,  , \, \Mc(\jb) = - \epsilon a- \delta b }  |j_1|^s |u_{j_1}(t)| \cdots |j_2|^s |u_{j_{i-2}}(t)| \dd t.
\end{multline*}
Multiplying by $(z_q - y_q)$ and summing in $k$ and $j$, we obtain 
$$
\Norm{Z_{P_i}(z) - X_{P_i}(y)}{\ell^1_s}
\leq 4 i (i-1)^s \Norm{P_i}{} \Big(\int_0^1 \Norm{u(t)}{\ell^1_s}^{i-2}  \dd t \Big) \Norm{z - y}{\ell^1_s}. 
$$
Hence we obtain the result after summing in $i$, using the fact that  $$\Norm{ty + (1-t)z}{\ell_s^1} \leq \max(\Norm{y}{\ell^1_s},\Norm{z}{\ell^1_s}).$$

Assume now that $P$ and $Q$ are homogeneous polynomials of degrees $k$ and $\ell$ respectively and with coefficients $a_\kb$, $\kb \in \Ic_k$ and $b_\ellb$, $\ellb \in \Ic_\ell$. It is clear that $\{P,Q\}$ is a monomial of degree $k + \ell - 2$ satisfying the zero momentum condition. Furthermore writing
$$
\{P,Q\}(z) = \sum_{ \jb \in \Ic_{k+\ell-2} } c_{\jb} z_\jb,
$$
$c_\jb$ expresses  as a sum of coefficients $a_{\kb}b_\ellb$ for which there exists an $a \in \Nc$ and $\epsilon \in \{\pm 1\}$ such that 
$$
(a,\epsilon) \subset \kb \in \Ic_k \quad \mbox{and} \quad (a,-\epsilon) \subset \ellb \in \Ic_\ell, 
$$
and such that if for instance $(a,\epsilon) = k_1$ and $(a,-\epsilon) = \ell_1$, we necessarily have  $(k_2,\ldots,k_k,\ell_2,\ldots,\ell_\ell) = \jb$. Hence for a given $\jb$, the zero momentum  condition on $\kb$ and on $\ellb$ determines  the value of $\epsilon a$ which in turn determines the value of $(\epsilon,a)$ when $\Nc = \N^d$  and determines two possible value of $(\epsilon,a)$ when $\Nc = \Z^d$. \\
This proves \eqref{Ebrack} for monomials. If 
$$
P = \sum_{i = 2}^k P_j \quad \mbox{and}\quad Q = \sum_{j = 2}^\ell Q_j
$$
where $P_i$ and $Q_j$ are homogeneous polynomials of degree $i$ and $j$ respectively, then we have 
$$
P = \sum_{n = 2}^{k+\ell -2} \sum_{i + j - 2 = n} \{P_i,Q_j\}. 
$$
Hence by definition of $\Norm{P}{}$ (see \eqref{EnormP}) and the fact that all the polynomials $\{P_i,Q_j\}$ in the sum are homogeneous or degree $i + j - 2$, we have by the previous calculations
\begin{multline*}
\Norm{P}{} = \sum_{n = 2}^{k+\ell -2} \Norm{\sum_{i + j - 2 = n} \{P_i,Q_j\}}{} \leq 2 \sum_{n = 2}^{k+\ell -2} \sum_{i + j - 2 = n} ij \Norm{P_i}{}\Norm{Q_j}{}
\\[2ex]
\leq 2 k \ell \Big( \sum_{i = 2}^k \Norm{P_i}{} \Big) \Big( \sum_{j = 2}^\ell \Norm{Q_j}{} \Big) = 2 k \ell \Norm{P}{}\Norm{Q}{}
\end{multline*}
where we used  \eqref{Epotpot} for the last equality. 

The last assertion, as well as the fact that the Poisson bracket of two real Hamiltonian is real,  follow immediately from the definition of the Poisson bracket. 
\end{Proof}

\begin{remark}\label{Rk-zero}
The zero momentum condition was crucially used to prove \eqref{Echamp} and \eqref{Ebrack}. For instance instead of \eqref{takecare}, we would have without this condition\footnote{Here we assumed that $|a_\jb|$ does not depend on permutation on the multi index $\jb$.}  in the case where $s = 0$,
$$\Norm{X_{P_i}(z) }{\ell^1} \leq 2 i \sum_{\ell\in \Zc}\sum_{\jb \in \Zc^{i-1}}|a_{\jb\ell}| | z_{j_1} \cdots z_{j_{i-1}}| $$
and this last expression cannot be controlled by the $\ell^1$ norm of $z$ without an extra decreasing property on the $|a_{\jb\ell}|$. For instance we can  assume that $\sum_\ell  |a_{\jb\ell}|$ is uniformly bounded with respect to $\jb$. See also \cite{Greb07} in the case where the phase space is $\ell^2_s$ instead of $\ell^1_s$. \end{remark}

With the previous notations, we consider in the following Hamiltonian functions of the form
\begin{equation}
\label{Edecomp}
H(z) = H_0(z) + P(z) = \sum_{a \in \Nc} \omega_a I_a(z)+ P(z),
\end{equation}
where for all $a\in \Nc$, 
$$
I_a(z) = \xi_a \eta_a
$$
are the {\em actions}   and
 $\omega_a \in \R$ are  the associated frequencies. We assume
\begin{equation}
\label{Eboundomega}
\forall\, a \in \Nc, \quad |\omega_a| \leq C |a|^m
\end{equation}
for some constants $C > 0$ and $m > 0$. 
The Hamiltonian system \eqref{Eham2} then reads
\begin{equation}
\label{Eham3}
\left\{
\begin{array}{rcll}
\dot\xi_a &=& \displaystyle - i \omega_a \xi_a - i \frac{\partial P}{\partial \eta_a}(\xi,\eta), & a \in \Nc,\\[2ex]
\dot\eta_a &=& \displaystyle i \omega_a \eta_a + i \frac{\partial P}{\partial \xi_a}(\xi,\eta),& a \in \Nc.  
\end{array}
\right.
\end{equation}

\subsection{Examples}\label{sec:ex}

In this subsection, we present two examples of equations that can be put under the previous form. We mention that many other systems can be put under the previous form. In particular, we stress out that we do not need any non resonance assumption on the frequencies $\omega_a$ as required in an perturbative approach (see for instance \cite{Greb07,BG06}). 

\subsubsection{Nonlinear Schr\"odinger equation}

We first consider  non linear Schr\"odinger equations of the form
\begin{equation}
\label{E1}
i \partial_t \psi = - \Delta \psi +  \partial_2g(\psi,\bar \psi),\quad x \in \T^d
\end{equation}
where $g:\C^2 \to \C$ is a polynomial of order $r_0$. We assume that $g(z,\bar z) \in \R$, and that $g(z,\bar z) = \mathcal{O}(|z|^3)$. 
The corresponding Hamiltonian functional is given by 
$$
H(\psi,\bar\psi) = \int_{\T^d} \left( | \nabla \psi | ^2 + g(\psi,\bar\psi) \right) \, \dd x.
$$

Let $\phi_{a}(x) = e^{i a\cdot x}$, $a \in \Z^d$ be the Fourier basis on $L^2(\T^d)$. With the notation 
$$
\psi = \Big(\frac{1}{2\pi}\Big)^{d/2}\sum_{a\in \Z^d} \xi_{a} \phi_{a}(x) \quad \mbox{and}\quad
\bar \psi = \Big(\frac{1}{2\pi}\Big)^{d/2}\sum_{a\in \Z^d} \eta_{a}\bar \phi_{a}(x)\,,
$$
 the Hamiltonian associated with the equation \eqref{E1} can (formally) be written
\begin{equation}
\label{E2}
H(\xi,\eta) = 
\sum_{a \in \Z^d} \omega_a \xi_a \eta_a  + \sum_{r = 3}^{r_0}\, 
\sum_{\ab,\bb } P_{\ab\bb}\,  \xi_{a_1} \cdots \xi_{a_p}\eta_{b_1}\cdots \eta_{b_q}.
\end{equation}
Here  $\omega_{a} =|a|^2$, satisfying \eqref{Eboundomega} with $m = 2$, are the eigenvalues of the Laplace operator $-\Delta.$ As previously seen in Examples \ref{ex}, \ref{ex2}, \ref{ex3}, the nonlinearity $P=\int_{\T^d}g(\psi(x),\phi(x))\dd x$ is real, satisfies the zero momentum condition and belongs to $\Hc_s$ (as $g$ is polynomial).

In this situation, working in the space $\ell^1$ for $\xi$ corresponds to working in a subspace of bounded functions 
$\psi(x)$. Similarly the control of the $\ell^1_s$ norm of $\xi$ for $s \geq 0$ leads to a control of $\Norm{\nabla^s\psi}{L^\infty}$. 

\subsubsection{Nonlinear wave equation}

As a second  concrete example 
we consider a 1-d nonlinear wave equation
\begin{eqnarray}
\label{nlw}
u_{tt}-u_{xx} = g(u)\ ,\quad x\in (0,\pi)  \ , \ t\in \R \ ,
\end{eqnarray}
with Dirichlet boundary condition: $u(0,t)=u(\pi,t)=0$ for any $t$.
We assume that $g:\R \to \R$ is polynomial of order $r_0 - 1$ with a zero of order two at $u=0$.
Defining $v=u_{t}$, \eqref{nlw} reads
$$
\partial_{t}
\begin{pmatrix}
    u \\ v
\end{pmatrix}=
\begin{pmatrix}
    v \\ u_{xx}+g(u)
\end{pmatrix}.
$$
Furthermore, let $H: H^{1}(0,\pi)\times L^2(0,\pi)\mapsto \R$ be defined by
\begin{equation} \label{Hnlw}
H(u,v)=\int_{S^1}\left( \frac 1 2 v^2 + \frac 1 2 u_{x}^2 +
G(u)\right) \dd x
\end{equation}
where $G$ such that
$\partial_{u}G=-g$ is a polynomial of degree $r_0$, then \eqref{nlw} expresses as a Hamiltonian system
\begin{eqnarray} \nonumber
    \partial_{t}
    \begin{pmatrix}
	u \\ v
    \end{pmatrix}&=&
    \begin{pmatrix}
	0&1 \\ -1 &0
    \end{pmatrix}
    \begin{pmatrix}
    -u_{xx}+ \partial_{u}G\\ v
    \end{pmatrix}\\[2ex] \label{nlwh}
    &=& J \nabla_{u,v}H(u,v)
    \end{eqnarray}
where $J=\begin{pmatrix}
	0&1 \\ -1 &0
    \end{pmatrix}$ represents the symplectic structure and where
    $\nabla_{u,v}=\begin{pmatrix}
    \nabla_{u} \\ \nabla_{v}
    \end{pmatrix}$ with
    $\nabla _u$ and $\nabla _v$ denoting the $L^2$ gradient with
    respect to  $u$ and  $v$  respectively. \\
Let $-\Delta_D$ be the Laplace operator with Dirichlet boundary conditions.  Let $  A= (-\Delta_D)^{1/2}$. We introduce the variables $  
(p,q)$ given by
$$
q:=A^{1/2}u\quad\mbox{and}\quad p:= A^{-1/2}v.
$$
Then, on $H^s(0,\pi)\times H^s(0,\pi)$ 
with $s\geq 1/2$, the Hamiltonian \eqref{Hnlw} takes the form $H_0+P$ with 
\begin{equation}
\label{Enlw1}
H_0(q,p)=\frac{1}{2} \big(\langle Ap,p \rangle_{L^2} +\langle
Aq,q \rangle_{L^2} \big)
\end{equation}
and 
\begin{equation}
\label{Enlw2}
P(q,p)=\int_{S^1}G(A^{-1/2}q) \dd x.
\end{equation}
In this context $\Nc=\N\setminus \{0\} $, 
$\omega_a=a$, ${a\in\Nc}$ are the eigenvalues of $A$  and $\phi_a=\sin ax$, $a\in\Nc$, the
associated eigenfunctions.

Plugging the decompositions 
$$
q(x) = \sum_{a \in \Nc} q_a \phi_a(x) \quad \mbox{and} \quad p(x) = \sum_{a\in \Nc} p_a \phi_a(x)
$$
into the Hamiltonian functional, we see that it takes the form 
$$
H = \sum_{a \in \Nc}\omega_a \frac{p_a^2 + q_a^2}{2} + P
$$
where $P$ is a function of the variables $q_a$. 
Using the complex coordinates 
$$
\xi_a = \frac{1}{\sqrt{2}} ( q_a + ip_a) \quad\mbox{and}\quad 
\eta_a = \frac{1}{\sqrt{2}} ( q_a - ip_a)
$$
the Hamiltonian function can be written under the form \eqref{E2} with a nonlinearity depending on $G$. In this case, the space $\ell^1$ for $z=(\xi,\eta)$ corresponds to functions $u(x)$ such that the Fourier transform $\widehat{u}(a) = \frac1\pi \int_{0}^\pi u(x) \sin(ax) \dd x$ satisfies $(a \widehat{u})_{a \in \Nc} \in \ell^1(\Nc)$ \quad \mbox{and}\quad $(a^{-1} \widehat{u})_{a \in \Nc} \in \ell^1(\Nc)$. This implies in particular a control of $u(x)$ and $\partial_x u(x)$  in $L^\infty(0,\pi)$.   More generally,  with a $z$ in some $\ell^1_s$ space, $s \in \N$, is associated a function $u(x)$ such that $\partial_x^k u(x) \in L^\infty(0,\pi)$ for $k = 0,\ldots, s+1$.

\subsection{Splitting schemes}
In this subsection we describe some Splitting schemes to which we will apply our technic in the next sections. The standard Lie-Trotter splitting methods for the PDE associated with the Hamiltonian $H_0+P$ consists in replacing the flow generated by $H$ during the time $h$ (the small time step) by the composition of the flows generated by $H_0$ and $P$ during the same time, namely
$$
\Phi_{H_0}^h \circ \Phi^h_P \quad \mbox{and} \quad  \Phi^h_P\circ \Phi_{H_0}^h.
$$
As explained in the introduction, it turns out that it is convenient to consider more general splitting methods including  in particular a regularization in the high modes of the linear part. Thus we replace the operator $hH_0$ by a more general Hamiltonian $A_0$.
Precisely let $\alpha_h(x)$ be a real function, depending on the stepsize $h$, satisfying $\alpha_h(0) = 0$ and $\alpha_h(x) \simeq x$ for small $x$. 
We define the diagonal operator $A_0$ by the relation
\begin{equation}
\label{Ealphah}
\forall\, j = (a,\delta) \in \Zc, \quad  A_0 z_j = \delta\alpha_h(h\omega_a) z_j. 
\end{equation}
 For $a \in \Nc$, we set $\lambda_a = \alpha_h(h\omega_a)$. 
We consider the splitting methods
\begin{equation}
\label{Esplis}
\Phi_P^h  \circ  \Phi_{A_0}^1 \quad\mbox{and}\quad \Phi_{A_0}^1\circ\Phi_P^h  
\end{equation}
where $\Phi_P^h$ is the exact flow associated with the Hamiltonian $P$, and where $\Phi_{A_0}^1$ is defined by the relation 
$$
\forall\, j = (a,\delta) \in \Zc,\quad 
\big(\Phi_{A_0}^1(z)\big)_j = \exp(i \delta \lambda_a) z_j.
$$
 We will mainly consider the cases listed in the table 1 below. 

\begin{table}[ht]
\renewcommand{\arraystretch}{1.9}
\begin{center} 
\begin{tabular}{|@{$\quad$}c@{$\quad$}|@{$\quad$}c@{$\quad$}|}
\hline
\textbf{\textsl{Method}} & \textbf{\textsl{$\alpha_h(x)$}} \\
\hline 
\sl Splitting & $\alpha_h(x) = x$  \\
\hline 
\sl Splitting + CFL & $\alpha_h(x) = x \mathds{1}_{x < c}(x)$  \\ 
 
\hline
\sl Mid-split & $\alpha_h(x) = 2\arctan(x/2)$\\ 
\hline
\sl Mid-split + CFL & $\alpha_h(x) = 2\arctan(x/2)\mathds{1}_{x < c}(x)$\\ 
\hline
\sl New scheme (I) & $\alpha_h(x) = h^{\beta}\arctan(h^{- \beta}x)$ \\
\hline 
\sl New scheme (II) & $\displaystyle\alpha_h(x) = \frac{x + x^2/h^\beta}{1 + x/h^\beta + x^2/h^{2\beta}}$ \\
\hline
\end{tabular}
\caption{Splitting schemes}
\label{Table1}
\end{center}
\end{table}

Let us comment these choices. The ``mid-split'' cases correspond to the approximation of the system
$$
\left\{
\begin{array}{rcll}
\dot\xi_a &=& \displaystyle - i \omega_a \xi_a & a \in \Nc,\\[2ex]
\dot\eta_a &=& \displaystyle i \omega_a \eta_a & a \in \Nc,
\end{array}
\right.
$$
by the midpoint rule (see also \cite{Ascher,Stern}). Starting from a given point $(\xi_a^0,\eta_a^0)$ we have by definition for the  first equation
$$
\xi^{1}_a = \Big(\frac{1 - ih\omega_a/2}{1 + ih\omega_a/2}\Big) \xi^{0}_a = \exp(-2i \arctan(h\omega_a/2) )\xi_a^0
$$
which is the solution at time $1$ of the system
$$
\left\{
\begin{array}{rcll}
\dot\xi_a &=& \displaystyle - i 2 \arctan(h\omega_a/2) \xi_a & a \in \Nc,\\[2ex]
\dot\eta_a &=& \displaystyle i 2\arctan(h\omega_a/2) \eta_a & a \in \Nc.  
\end{array}
\right.
$$
Thus in this case the Hamiltonian $A_0$  is given by 
$$
A_0(\xi,\eta) = \sum_{a\in \Zc} 2 \arctan(h\omega_a/2) \, \xi_a\eta_a. 
$$
Note that using the relation 
\begin{equation}
\label{Eestatan}
\forall\, y \in \R, \quad 
| \arctan(y)- y | \leq \frac{|y|^3}{3}
\end{equation}
we obtain for all $a \in \Nc$,
$$
 |\exp(-i h \omega_a) -  \exp(-2i\arctan(h\omega_a/2)) | \leq C h^3 \omega_a^3
$$
for some constant $C$ independent of $a$. 
Using the bound \eqref{Eboundomega}, we get for all $z$, 
$$
\Norm{\Phi_{H_0}^h(z) - \Phi_{A_0}^1(z)}{L^2} \leq C h^3 \Norm{z}{H^{3m}}. 
$$
More generally, we have the following approximation result: 

\begin{lemma}
\label{Pordre}
Assume that the function $\alpha_h(x)$ satisfies: 
\begin{equation}
\label{Econdalph}
\forall\, x > 0, \quad  |\alpha_h(x) -x | \leq C h^{-\sigma} x^{\gamma}
\end{equation}
for some constants $C> 0$, $\sigma \geq 0$ and $\gamma \geq 2$. Then we have 
\begin{equation}
\label{EapproxA}
\Norm{\Phi_{H_0}^h(z) - \Phi_{A_0}^1(z)}{L^2} \leq C  h^{\gamma - \sigma} \Norm{z}{H^{\gamma m}}. 
\end{equation}
\end{lemma}
\begin{Proof}
For a given $a \in \Nc$, we have 
$$
|\alpha_h(h\omega) -h\omega_a | \leq C h^{ \gamma -\sigma} \omega_a^{\gamma}.
$$
Hence owing to the fact that $|e^{ix} - e^{iy}| \leq |x - y|$ for real $x$ and $y$, we obtain
$$
|\exp(-i h \omega_a) - \exp(-i \alpha_h(h\omega_a) |
\leq  C  h^{ \gamma -\sigma} \omega_a^{\gamma},
$$
and this yields the result. 
\end{Proof}


\begin{commentary}
Under the assumption that $P$ acts on sufficiently high index Sobolev spaces $H^s$, the previous result can be combined with standard convergence analysis to show that the splitting methods \eqref{Esplis} yield consistent approximation of the exact solution $\Phi_H^h$ provided the initial solution is smooth enough (depending on $m$). The condition $\gamma - \sigma \geq 2$ guarantees a local order $2$ in \eqref{EapproxA} that will be of the same order of the error made by the splitting decomposition after one step. Such a local error propagates to a global error of order $1$, which means that for a give finite time $T$, the error with the exact solution after $n$ iterations with $nh = T$ will be of order $n\times h^2 \simeq h$ up to constants depending on $T$, and under the assumption that the numerical solution remains smooth. Give more precise results would be out of the scope of this paper, and we refer to \cite{L08a} for the case of NLS. 
\end{commentary}

Let us consider 
the function 
\begin{equation}
\label{NS1}
\alpha_h(x) = h^{\beta}\arctan(h^{ - \beta}x)
\end{equation}
for $1 > \beta \geq 0$. It satisfies \eqref{Econdalph} with $C = 1/3$, $\sigma = 2 \beta$ and $\gamma = 3$ (see \eqref{Eestatan}). 
Hence for $\beta = 1/2$, the estimate \eqref{EapproxA} shows a local error of order $\gamma - \sigma = 2$, and hence the splitting schemes \eqref{Esplis} remains of local order $2$ (though with more smoothness required than with the midpoint approximation) which means that the error made after one step is of order $h^2$. Of course, when $\beta = 1$, the approximation is not consistent (local error of order $1$). 

The second example 
\begin{equation}
\label{NS2}
\alpha_h(x) = \frac{x + x^2/h^\beta}{1 + x/h^\beta + x^2/h^{2\beta}}
\end{equation}
exhibits similar properties. Note that the simple choice
\begin{equation}
\label{NS3}
\alpha_h(x) = \frac{x}{1 + x/h^\beta}
\end{equation}
ensures only a local error of order $h^{2 - \beta}$ ($\gamma = 2$ and $\sigma = \beta$ in \eqref{EapproxA}). Hence the corresponding splitting schemes \eqref{Esplis} are of global order $1 - \beta$ (hence $1/2$ in the case where $\beta = 1/2$). 

All these ``new'' schemes have the particularity that $\alpha_h(x) \simeq x$ when $x$ is small, but when $x \to \infty$, we have $\alpha_h(x) \simeq h^{\beta}$. Their use thus leads to a stronger regularization effect in the high modes than the midpoint approximation, without breaking the order of approximation for smooth functions. We will see in section \ref{formal} and \ref{analytic} that this property  allows to construct a modified equation over exponentially long time for all these schemes. 

Note that in practice, the implementation of the schemes associated with the filter functions \eqref{NS1} or \eqref{NS2} {\em a priori} requires the knowledge of a spectral decomposition of $H_0$. This will be the case for NLS on the torus, or NLW with Dirichlet boundary conditions, where the switch from the $x$-space (to calculate $\Phi^h_P$) to the Fourier space (to calculate $\Phi_{A_0}^1$) can be easily implemented using fast Fourier transformations. 


%

\section{Recursive equations}\label{formal}

In this section we explain the strategy in order to prove the existence of a modified energy. We will see that it leads to solve by induction a sort of   homological equation in the spirit of normal form theory (see for instance \cite{Greb07}). For simplicity, we consider only the splitting method $\Phi_P^h \circ \Phi_{A_0}^1$. The second Lie splitting  $\Phi_{A_0}^1\circ \Phi_P^h $ can be treated similarly. 

We look for a real Hamiltonian function $Z(t,\xi,\eta)$ such that for all $t \leq h$ we have
\begin{equation}
\label{Ehhh}
\Phi_P^t \circ \Phi_{A_0}^1 = \Phi_{Z(t)}^1
\end{equation}
and such that $Z(0) = A_0$. 

For a given Hamiltonian $K \in \Hc_s$, we denote by $\Lc_K$ is the Lie differential operator associated with the Hamiltonian vector field $X_K$: for a given function $g$ acting on  $\ell^1_s$, $s \geq 0$,  and taking values on $\C$ or $\ell^1_s$, we have
$$
\Lc_K(g) = \sum_{j \in \Zc} (X_K)_j \frac{\partial g}{\partial z_j}. 
$$
Denoting by $z(t)$ the flow generated by $X_K$ starting from $z\in \ell^1_s$, i.e.\ $z(t)=\Phi^t_K(z)$, we have (if $K$ is in  $\Hc_s$) 
$$
z^{(k)}(t)= \Lc_K^k[I](z(t)) \quad \mbox{for all}\quad k\in \N,
$$
where $I$ define the identity vector field: $I(z)_j=z_j$. 
Thus we can write at least formally
\begin{equation}
\label{exp}
\Phi_{K}^1 = \exp( \Lc_K)[I]. 
\end{equation}
Differentiating the exponential map we calculate as in  \cite[Section III.4.1]{HLW}
$$
\frac{\dd}{\dd t} \Phi_{Z(t)}^1 = X_{Q(t)} \circ \Phi_{Z(t)}^1,
$$
where the differential operator associated with $Q(t)$ is given by 
$$
\Lc_{Q(t)} = \sum_{k \geq 0} \frac{1}{(k+1)!} \mathrm{Ad}^k_{\Lc_{Z(t)}} (\Lc_{Z't)})
$$
with
$$
\mathrm{Ad}_{\Lc_A}(\Lc_H) = [\Lc_A,\Lc_H]
$$ 
the commutator of two vector fields. 

As the vector fields are Hamiltonian, we have
$$
[\Lc_A,\Lc_H] = \Lc_{\{A,H\}}, 
$$
 where 
$$
\mathrm{ad}_{K}(G) = \{K,G\}. 
$$
Hence we obtain
the formal series equation for $Q$:
\begin{equation}
\label{EQt}
Q(t) = \sum_{k \geq 0                                                                                                                                                                                                                                                                                                                                                                                                                                                                                                                                                                                                                                                                                                                                                                                                                                                                                                                                                                                                                                                                                                                                                                                                                                                                                                                                                                                                                                                                                                                                                                                                                                                                                                                                                                                                                                                                                                                                                                                                                                                                                                                                                                                                                                                                                                                                                                                                                                                                                                                                                                                                                                                                                                                                                                                                                                                                                                                                                                                                                                                                                                                                                                                                                                                                                                                                                                                                                                                                                                                                                                                                                                                                                                                                                                                                                                                                                                                                                                                                                                                                                                                                                                                                                                                                                                                                                                                                                                                                                                      } \frac{1}{(k+1)!} \mathrm{ad}^k_{Z(t)} Z'(t)
\end{equation}
where $Z'(t)$ denotes the derivative with respect to $t$ of the Hamiltonian function $Z(t)$. 

Therefore taking the derivative of \eqref{Ehhh}, we obtain 
$$
X_P \circ \Phi_P^t \circ \Phi_{A_0}^1  = X_{Q(t)} \circ \Phi_{Z(t)}^1
$$
and hence the equation to be satisfied by $Z(t)$ reads: 
\begin{equation}
\label{Ehomo}
\sum_{k \geq 0} \frac{1}{(k+1)!} \mathrm{ad}^k_{Z(t)} Z'(t) = P. 
\end{equation}
Notice that the series $\sum_{k\geq 0} \frac{1}{(k+1)!} z^k= \frac{e^z-1}{z}$ is invertible in the open disc $|z|<2\pi$ with inverse given by $\sum_{k \geq 0} \frac{B_k}{k!}z^k$ where $B_k$ are the Bernoulli numbers. 
So formally,  Equation \eqref{Ehomo} is equivalent to the formal series equation (see also \cite{DF09}, Eqn.\ (3.1))
\begin{equation}\label{Ehomo2}
Z'(t) = \sum_{k \geq 0} \frac{B_k}{k!} \mathrm{ad}^k_{Z(t)} P. 
\end{equation}

Plugging an Ansatz expansion $Z(t) = \sum_{\ell \geq 0} t^\ell Z_\ell$ into this equation, we get $Z_0 = A_0$ and for $n \geq 0$
\begin{equation}
\label{Erec}
(n+1) Z_{n+1} =  \sum_{k \geq 0 }\frac{B_k}{k!} \sum_{\ell_1 + \cdots + \ell_k = n} \mathrm{ad}_{Z_{\ell_1}} \cdots \mathrm{ad}_{Z_{\ell_k}} P. 
\end{equation}

\begin{commentary}
The analysis made to obtain this recursive equation is formal. To obtain our main result, we will  verify that the series we manipulate are in fact  convergent series in $\Hc_s$ uniformly on balls of $\ell^1_s$ that contains the different flows  involved in the formulas (see in particular Lemma \ref{Linvert} below).\\
For instance \eqref{exp} holds true as soon as $z(t)$ remains in a ball  $B^s_M:=\{z\in \ell^1_s \, | \,  \Norm{z}{\ell^1_s}\leq M\}$ for $0\leq t\leq 1$, and that the series $\sum\frac{\Lc^k_K[I](z)}{k!}$ is uniformly convergent on $B^s_M$. Notice that  this in turn implies that $t\mapsto z(t)$ is analytic on the complex disc of radius $1$.
\end{commentary}
\begin{commentary}
In the case of Strang splitting methods of the form
$$
\Phi_P^{h/2} \circ \Phi_{A_0}^1 \circ \Phi_P^{h/2}
$$
we can apply the same strategy and look for a Hamiltonian $Z(t)$ satisfying for all $t \leq h$, 
$$
\Phi_{Z(t)}^1 = \Phi_P^{t/2} \circ \Phi_{A_0}^1 \circ \Phi_P^{t/2}
$$
and this yields to equations similar to \eqref{Erec}. We do not give the details here, as the analysis will be very similar. 
\end{commentary}
The key lemma in order to prove that the previous series converge (and thus to justify the previous formal analysis) is the following one (and whose proof is straightforward calculus):

\begin{lemma}
\label{L1}
Assume that 
$$
Q(z) = \sum_{j \in \Zc} a_{\jb} z_\jb
$$
is a polynomial, then 
$$
\mathrm{ad}_{Z_0}(Q) = \sum_{j \in \Zc} i\Lambda(\jb)a_\jb z_\jb
$$
where for a multi-index $\jb = (j_1,\ldots, j_r)$ with for $i = 1,\ldots,r$, $j_i = (a_i,\delta_i) \in \Nc\times\{\pm 1\}$, we set 
$$
\Lambda(\jb) = \delta_1 \lambda_{a_1} + \cdots + \delta_r \lambda_{a_r}. 
$$
\end{lemma}
Hence we see that if $\Lambda(\jb) < 2\pi$ we will be able to define at least the first term $Z_1$ by summing the series in $k$  in the formula \eqref{Erec}. 

\section{Analytic estimates}\label{analytic}

We assume in this section that $\alpha_h$ and $h$ satisfy the following condition: there exist $r$ and a constant $\delta < 2\pi$ such that 
\begin{equation}
\label{Enonres}
\forall\, n \leq r\, \quad \forall\, \jb \in \Ic_n,\quad 
|\Lambda(\jb)|\leq 2\pi - \delta. 
\end{equation}
Using Lemma \ref{L1}, this condition implies that for any polynomial $Q \in \Pc_r$,  we have the estimate 
\begin{equation}
\label{Enonres3}
\Norm{\ad_{Z_0} Q}{} \leq (2\pi - \delta) \Norm{Q}{}
\end{equation}
as for homogeneous polynomials, the degree of $\ad_{Z_0}Q$ is the same as the degree of $Q$. 

\begin{theorem}\label{T1}
Let $r_0 \geq 3$. Assume that $P \in \Pc_{r_0}$ and that the condition \eqref{Enonres} is fulfilled for some constants $\delta$ and $r$. Then for $n \leq N:=\frac{r-2}{r_0 - 2}$ we can define polynomials $Z_n \in \Pc_{n(r_0 - 2)+2}$ satisfying the equations \eqref{Erec} up to the order $n$, and satisfying the estimates $\Norm{Z_1}{} \leq c$ and for $2\leq n \leq N$, 
\begin{equation}
\label{EestZn}
\Norm{Z_n}{} \leq c(C n)^{n-2} 
\end{equation}
for some constants $c$ and $C$ depending only on $\Norm{P}{}$, $r_0$ and $\delta$. If moreover 
$P \in \Sc\Pc_{r_0}$ then $Z_n \in \Sc\Pc_{n(r_0-2)+2}.$
\end{theorem}

\begin{Proof}
Let
$$
P(z) = \sum_{\ell = 2}^{r_0} \sum_{\jb \in \Ic_\ell} a_{\jb} z_{\jb}.
$$
First we prove the existence of the $Z_k$ for $k\leq N$. The equation \eqref{Erec} for $n = 0$ reads
$$
Z_1 =   \sum_{k \geq 0 }\frac{B_k}{k!}  \mathrm{ad}_{Z_{0}}^k P. 
$$
The previous Lemma and the condition \eqref{Enonres} show that $Z_1$ exits and is given by
$$
Z_1 = \sum_{\ell = 2}^{r_0} \sum_{\jb \in \Ic_\ell} \frac{i\Lambda(\jb)}{\exp(i \Lambda(\jb)) - 1} a_{\jb} z_{\jb}.
$$
Further we deduce immediately that $Z_1$ is real, and satisfies $\Norm{Z_1}{} \leq c_\delta \Norm{P}{}$ for some constant $c_\delta$. 

 Assume now that the $Z_k$ are constructed for $0 \leq k \leq n$, $n \geq 1$ and are such that $Z_k$ is a polynomial of degree $k(r_0 - 2)+2$. Formally $Z_{n+1}$ is  defined as a series
 $$
Z_{n+1}  = \frac{1}{n+1} \sum_{k \geq 0 }\frac{B_k}{k!} A_k
$$
where
$$
A_k=\sum_{\ell_1 + \cdots + \ell_k = n} \mathrm{ad}_{Z_{\ell_1}} \cdots \mathrm{ad}_{Z_{\ell_k}} P.
$$
Let us prove that this series converges absolutely. In the previous sum, we separate the number of indices $j$ for which $\ell_j = 0$. For them, we can use \eqref{Enonres3}. Only for  the other indices, we will use the estimates of Proposition \ref{P1} by taking into account that the right-hand side  is a sum of terms that are all real polynomials of degree $(\ell_1 + \cdots + \ell_k)(r_0 - 2) + r_0 = (n+1)(r_0 - 2) + 2 \leq (n+1)r_0$ and hence the inequality of Proposition \ref{P1} is only used with polynomials of order less than $(n+1)r_0$. Thus we write for $k\geq n$
\begin{multline*}
\Norm{A_k}{}:= \Norm{\sum_{\ell_1 + \cdots + \ell_k = n} \mathrm{ad}_{Z_{\ell_1}} \cdots \mathrm{ad}_{Z_{\ell_k}} P}{}\leq
 \\[2ex]
 \sum_{i = 1}^n \frac{k! \ (2\pi - \delta)^{k-i}}{(k-i)!\ i!}\sum_{\ell_1 + \cdots + \ell_i = n | \ell_j > 0} 
(n+1)^{i - 1 } 2^i r_0^{2i}\ell_1 \Norm{Z_{\ell_1}}{} \cdots \ell_i \Norm{Z_{\ell_i}}{}  \Norm{P}{}
 \\[2ex]
\leq  (2\pi - \delta)^{k-n}k^n\sum_{i = 1}^n \sum_{\ell_1 + \cdots + \ell_i = n | \ell_j > 0} 
(n+1)^{i - 1 } 2^i r_0^{2i}\ell_1 \Norm{Z_{\ell_1}}{} \cdots \ell_i \Norm{Z_{\ell_i}}{}  \Norm{P}{},
 \end{multline*}  
 and thus $\sum_{k \geq 0 }\frac{B_k}{k!} A_k$ converges and  $Z_{n+1}$ is well defined up to $n+1\leq N$.  
 
 Further, as the right-hand side of \eqref{Erec} is a sum of term that are all real polynomials of degree $(n+1)(r_0 - 2) + 2$, $Z_{n+1}$ is a polynomial of degree $(n+1)(r_0 - 2) + 2$.
 
Now we have to prove the estimate \eqref{Enonres3}. Following the previous calculation we get 
\begin{multline*}
(n+1)\Norm{ Z_{n+1}}{}  \leq 
\\[2ex]  \sum_{i = 1}^n \sum_{\ell_1 + \cdots + \ell_i = n | \ell_j > 0} \sum_{k \geq i }\frac{B_k}{i!\ (k-i)!} (2\pi - \delta)^{k-i}
(n+1)^{i - 1 } 2^i r_0^{2i} \times \\[2ex]\ell_1 \Norm{Z_{\ell_1}}{} \cdots \ell_i \Norm{Z_{\ell_i}}{}  \Norm{P}{}.
\end{multline*}
On the other hand, the entire series $f(z):=\sum_{k \geq 1 }\frac{B_k}{k!} z^k$ defines an analytic function on the disc $|z|< 2\pi$. Thus its $i-th$ derivative $\sum_{k \geq i }\frac{B_k}{(k-i)!} z^{k-i}$ also defines an analytic function on the same disc and, by Cauchy estimates, there exist a constant $C_\delta=\sup_{|z|\leq 2\pi-\delta/2}|f(x)|$ such that 
$$
\sum_{k \geq i }\frac{B_k}{(k-i)!} (2\pi - \delta)^{k-i} \leq C_\delta\ i!\ \left(\frac{2}{\delta}\right)^i.
$$
We then define for $n \geq 1$
$$
\zeta_n = nr_0  \Norm{Z_n}{}.
$$
These numbers satisfy the estimates, for $n \geq 0$, 
$$
\zeta_{n+1} \leq \delta_n^0 r_0 c_\delta\Norm{P}{} 
+ C_\delta  r_0 \Norm{P}{} \sum_{i = 1}^n \sum_{\ell_1 + \cdots + \ell_i = n | \ell_j > 0}  
(n+1)^{i -1 } (4\delta^{-1}r_0)^{i} \zeta_{\ell_1} \cdots \zeta_{\ell_i}. 
$$
Let us fix $N \geq 1$. We have for $n = 0,\ldots,N$ 
\begin{multline*}
4\delta^{-1}r_0(N+1)\zeta_{n+1} \leq 4\delta^{-1}r_0(N+1) \delta_n^0 r_0 c_\delta\Norm{P}{} \\[2ex]
+  C_\delta4\delta^{-1}r_0 ^2 \Norm{P}{} \sum_{i = 1}^n \sum_{\ell_1 + \cdots + \ell_i = n | \ell_j > 0} 
(N+1)^{i } (4\delta^{-1}r_0)^{i} \zeta_{\ell_1} \cdots \zeta_{\ell_i}.
\end{multline*}
Let $\beta_j$, $j = 0,\ldots,N$ the sequence satisfying
\begin{equation}
\label{Ebeta}
\beta_{n+1} =  \delta_n^0 C_1 + C_2 \sum_{i = 1}^n \sum_{\ell_1 + \cdots + \ell_i = n | \ell_j > 0}   \beta_{\ell_1} \cdots \beta_{\ell_i}
\end{equation}
where
$$
C_1 = 4(N+1) r_0^2c_\delta\ \delta^{-1}\Norm{P}{} \quad  \mbox{and}\quad C_2 = \frac{4C_\delta}{ \delta} r_0^2 \Norm{P}{}. 
$$ 
By induction, we have that for all $n = 0,\ldots,N$, 
$$
(N+1)\frac{4r_0}{\delta} \zeta_n \leq \beta_n. 
$$
Multiplying \eqref{Ebeta} by $t^{n+1}$ and summing in $n\geq0$ we see that the formal series $\beta(t) = \sum_{j \geq 1} t^j \beta_j$ satisfies the relation
$$
\beta(t) = tC_1 + tC_2 \left( \frac{1}{1 - \beta(t)} - 1\right).
$$
This yields 
$$
(1 - \beta(t))(\beta(t) - tC_1) = tC_2 \beta(t)
$$
or equivalently 
$$
\beta(t)^2 - \beta(t) (1 + t(C_1-C_2)) + t C_1=0. 
$$
The discriminant of this equation is
$$
(1 + t(C_1-C_2))^2 - 4 t C_1 = 1 - 2t (C_1 +  C_2) + t^2(C_1-C_2)^2
$$
and hence, for $t\leq 1/2(C_1+C_2)$, we find using $\beta(0) = 0$, 
$$
2\beta = 1 + t(C_1-C_2) - \Big(1 - 2t(C_1 +  C_2) + t^2(C_1-C_2)^2\Big)^{1/2}.
$$
We verify that for $t \leq 1/2(C_1 - C_2)$ we have 
$$
2\beta \leq   \frac 3 2. 
$$
By analytic estimate, we obtain that for all $n\geq 0$ we have
$$
\beta_{n} = \frac{\beta^{(n)}(0)}{n!}\leq \frac{3}{2} \Big(2(C_1 - C_2)\Big)^{n}.
$$
For $n = N$, this yields
$$
\beta_N \leq (CN \Norm{P}{})^{N}
$$
for some constant $C$ depending on $r_0$ and $\delta$. We deduce the claimed result from the expression of $\zeta_N$. 
\end{Proof}

 For $s \geq 0$, we define 
$$
B^s_M = \{\, z \in \ell^1_s \, | \, \Norm{z}{\ell^1_s} \leq M\,\}.
$$
and we will use the notation $B_M = B^0_M$. 
\begin{theorem}
\label{TB}
Let $r_0 \geq 3$, $s \geq 0$  and $M \geq 1$ be  fixed.  We assume that $P \in \Pc_{r_0}$ and that the condition \eqref{Enonres} is fulfilled for some constants $\delta$ and $r\geq r_0$ and we denote by  $N$ the largest integer smaller than 
$\frac{r-2}{r_0 - 2}. $
Then there exist  constants $c_0$ and $C$ depending   on $r_0$, $s$, $\delta$, $\Norm{P}{}$ and $M$   such that for all $hN \leq c_0$,  there exists a real Hamiltonian polynomial $H_h \in \Pc_{N(r_0 - 2)+2}$ such that for all $z \in B^s_M$, we have 
\begin{equation}\label{Estfond}
\Norm{\Phi_P^h \circ \Phi_{A_0}^1(z) - \Phi_{H_h}^h(z) }{\ell^1_s} \leq h^{N+1} (CN)^{N}. 
\end{equation}
Moreover, assuming that 
$$
P(z) = \sum_{\ell = 1}^{r_0} \sum_{\jb \in \Ic_\ell} a_{\jb} z_{\jb}
$$
then for $z \in B^s_M$ we have 
\begin{equation}
\label{EestH1}
|H_h(z) - H^{(1)}_h(z)| \leq Ch 
\end{equation}
where 
\begin{equation}
\label{EH1}
H^{(1)}_h(z) = \sum_{a \in \Nc} \frac{1}{h}\alpha_h(h\omega_a) \xi_a\eta_a + \sum_{\ell = 1}^{r_0} \sum_{\jb \in \Ic_\ell} \frac{i\Lambda(\jb)}{\exp(i \Lambda(\jb)) - 1} a_{\jb} z_{\jb}. 
\end{equation}
If finally, $P \in \Sc\Pc_{r_0}$, then $H_h \in \Sc\Pc_{N(r_0 - 2)+2}$. 
\end{theorem}

\begin{Proof}
We define the real Hamiltonian $H_h=\frac{Z_h(h)}{h}$, where
$$
Z_h(t) = \sum_{j = 0}^N t^j Z_j,
$$
 and where, for $j=0,\cdots,N$, the polynomials $Z_j$  are defined in Theorem \ref{T1}. Notice that  $N$ depends on $r$ and thus on $h$ via the condition \eqref{Enonres}.\\
By definition, $Z_h(t)(z)$ is a polynomial of order $N(r_0 - 2)+2\leq Nr_0$ and using Theorem \ref{T1} we get 
$$
\Norm{Z_h(t)}{} \leq c_1( 1 +  \sum_{j = 2}^N (C tj)^{j-2})  < \infty.
$$
for some constant $c_1$ depending on $\delta$ and $\Norm{P}{}$. 
Thus $Z_h \in \Pc_{N(r_0 - 2)+2}$ (and in $\Sc\Pc_{N(r_0 - 2)+2}$ if $P \in \Sc\Pc_{r_0}$).

We will use the fact that for all $s$, there exists a constant $c_s$ such that for all $j \geq 1$, 
\begin{equation}
\label{Ecs}
j^{s} \leq (c_s )^{j}. 
\end{equation}
Now, as for all $j$, $Z_j$ is a polynomial of order $j(r_0 - 2) + 2 \leq j{r_0}$ with a zero of order at least 2 in the origin, we have using Proposition \ref{P1} and Theorem \ref{T1} that  for $z \in B^s_{9M}$ and $j \geq 1$.
$$
\Norm{X_{Z_j}(z)}{\ell^1} \leq 2c(j r_0)^{s+1} (C{}j)^{j-1}\big(\sup_{k = 1,\ldots,jr_0 -1}\Norm{z}{\ell^1}^{k}\big)
\leq  M (2 C c_ s r_0^{s+1} c j  (9M)^{r_0})^j
$$
where the constants $c$ and $C$ are given by estimate \eqref{EestZn}. 
On the other hand we have using Lemma \ref{L1} and \eqref{Enonres} 
$$
\Norm{X_{Z_0}(z)}{\ell^1} \leq 2\pi \Norm{z}{\ell^1} \leq 2\pi M. 
$$
Hence, for $t  \leq (4NC c_s r_0^{s+1}c )^{-1}(9M)^{-r_0}$ we have 
\begin{equation}
\label{Egrad}
\Norm{X_{Z_h(t)}(z)}{\ell^1} \leq 2\pi M+  M   \sum_{j = 1}^N   (t 2 C c_s r_0^{s+1} cN  (9M)^{r_0})^j  < (2\pi +1) M<8M.
\end{equation}
Therefore by a classical bootstrap argument, the time 1  flow
$\Phi_{Z_h(t)}^1$ map $B^s_{M}$ into $B^s_{9M}$ provided that $t  \leq (4NC c_s r_0^{s+1}c)^{-1}(9M)^{-r_0}$.

On the other hand,  $\Phi^1_{A_0}$ is an isometry of $\ell^1_s$ and hence maps $B^s_M$ into itself, while using again Proposition \ref{P1}, we see that $\Phi^t_P$ maps $B^s_M$ into $B^s_{9M}$ as long  as $t\leq (4r_0^{s+1}\Norm{P}{}M^{(r_0-1)})^{-1}$. We then define
\begin{multline}
\label{EestT}
T\equiv T(N,M,r_0,s,\delta,\Norm{P}{})\\
:=\min\{ (4r_0^{s+1}\Norm{P}{}M^{(r_0-1)})^{-1}, (4NC c_s r_0^{s+1}c )^{-1}(9M)^{-r_0}\}
\end{multline}
and we assume in the sequel that $0\leq t\leq T$ in such a way that all the flows remain in the ball $B_{9M}$.\\
Let $u(t) = \Phi_P^t \circ \Phi_{A_0}^1 (z) - \Phi_{Z_h(t)}^1(z)$ and denote by $Q_h(t)$ the Hamiltonian defined by
$$
Q_h(t)=\sum_{k \geq 0} \frac{1}{(k+1)!} \mathrm{ad}^k_{Z_h(t)} Z'_h(t).
$$
\begin{lemma}
\label{Linvert}
For $t\leq T$ given in \eqref{EestT}, 
the Hamiltonian $Q_h(t) \in \Hc_s$ and satisfies  for $z \in B^s_M$
\begin{equation}
\label{Erelat}
\frac{\dd}{\dd t}\Phi_{Z_h(t)}^1(z)= X_{Q_h(t)}\circ \Phi^1_{Z_h(t)}(z).
 \end{equation}
\end{lemma}
We postpone the proof of this Lemma to the end of this section. 

Using this result, as $u(0) = 0$,  we get for 
$t\leq T$ given in \eqref{EestT}
\begin{align*}
\Norm{u(t)}{\ell^1_s} &\leq \int_0^t \Norm{X_P(\Phi_P^s \circ \Phi_{A_0}^1 (z))-X_{Q_h(s)}(\Phi_{Z_h(s)}^1(z))}{\ell^1_s}\dd s \\[2ex]
&\leq \int_0^t \Norm{X_P(\Phi_{Z_h(s)}^1(z))-X_{Q_h(s)}(\Phi_{Z_h(s)}^1(z))}{\ell^1_s}\dd s \\[2ex]
&+\int_0^t \Norm{X_P(\Phi_P^s \circ \Phi_{A_0}^1 (z))-X_P(\Phi_{Z_h(s)}^1(z))}{\ell^1_s}\dd s.
\end{align*}
Therefore for 
$t\leq T$
\begin{equation}\label{Est2}
\Norm{u(t)}{\ell^1_s}\leq \int_{0}^t \sup_{z\in B_{9M} }\Norm{X_P(z)-X_{Q_h(s)}(z)}{\ell^1_s} \dd s +L_P\int_0^t \Norm{u(s)}{\ell^1_s} \dd s
\end{equation}
where  using Equation \eqref{Elip} in  Proposition \ref{P1}, we can take
$$
L_P = 4 r_0^{s+1}  \Norm{P}{}(9M)^{r_0 - 2}.
$$ 
So it remains to estimate
$ \sup_{z\in B_{9M} }\Norm{X_P(z)-X_{Q_h(t)}(z)}{\ell^1_s}$ for $z\in B_{9M}$ and $t\leq T$. 

Now by definition of $Q_h(t)$ and using Lemma \ref{Linvert} we have 
$$
Z_h'(t) = \sum_{k = 0}^\infty \frac{B_k}{k!} {\rm ad}^k_{Z_h(t)} Q_h(t)
$$
where the right hand side actually defines a convergent series by the argument used in the proof of Theorem \ref{T1}. By construction (cf. section 3), we have
$$
\sum_{k = 0}^\infty \frac{B_k}{k!} {\rm ad}^k_{Z_h(t)} (Q_h(t) - P) = \mathcal{O}(t^{N})
$$
in the sense of Hamiltonian in the space $\Hc_s$. Taking the inverse of the series, we see 
\begin{equation}\label{Est1}
 Q_h(t)- P=\sum_{n\geq N}K_n
 \end{equation}
where we have the explicit expressions
\begin{equation}
\label{EKn}
K_n = \sum_{\ell + m = n| m < N}(m+1)  \sum_{k \geq 0 }\frac{1}{(k+1)!} \sum_{\ell_1 + \cdots + \ell_k = \ell | \ell_j \leq N} \mathrm{ad}_{Z_{\ell_1}} \cdots \mathrm{ad}_{Z_{\ell_k}} Z_{m+1}.
\end{equation}
Estimates similar to the one in the proof of Theorem \ref{T1} lead to 
\begin{multline*}
\Norm{K_n}{} \leq  \sum_{\ell + m = n | m < N}(m+1)  \sum_{i = 0}^\ell 2^i r_0^{2i} (n+1)^i \sum_{k \geq i }\frac{(2\pi-\delta)^{(k-i)}}{i!\ (k-i)!} \times \\[2ex]
\sum_{\ell_1 + \cdots + \ell_i = \ell | 0 < \ell_j \leq N} \ell_1 \Norm{Z_{\ell_1}}{} \cdots \ell_{i-1} \Norm{Z_{\ell_{i-1}}}{} \ell_i \Norm{Z_{\ell_i}}{} \Norm{Z_{m+1}}{}
\end{multline*}
and hence after summing in $k$,
\begin{multline*}
\Norm{K_n}{} \leq C_1 \sum_{\ell + m = n | m < N}(m+1)  \sum_{i = 0}^\ell \frac{2^i r_0^{2i} (n+1)^i}{i! } \times \\[2ex]
\sum_{\ell_1 + \cdots + \ell_i = \ell | 0 < \ell_j \leq N} \ell_1 \Norm{Z_{\ell_1}}{} \cdots \ell_{i-1} \Norm{Z_{\ell_{i-1}}}{} \ell_i \Norm{Z_{\ell_i}}{} \Norm{Z_{m+1}}{}
\end{multline*}
for some constant $C_1$ depending on $\delta$. 
Using the estimates in Theorem \ref{T1}, we have for $\ell_j > 0$ and $\ell_j \leq N$, 
$$
\ell_j \Norm{Z_{\ell_j}}{} \leq c(C \ell_j)^{\ell_j - 1} \leq c(CN)^{\ell_j - 1}. 
$$
Using moreover $\Norm{Z_{m+1}}{} \leq c(C N)^{m}$, and 
as  the number of integer $\ell_1,\ldots,\ell_i$ stricly positive such that $\ell_1 + \cdots + \ell_i = \ell$ is bounded by $2^{2\ell}$, we obtain
\begin{multline*}
\Norm{K_n}{} \leq c\, C_1 \sum_{\ell + m = n\ m < N} (m+1) \sum_{i = 0}^\ell \frac{(2c)^i r_0^{2i} (n+1)^i}{i! } \times \\[2ex]
\sum_{\ell_1 + \cdots + \ell_i = \ell | 0 < \ell_j \leq N} (C N)^{\ell_1 + \cdots + \ell_i + m - i}\\[2ex]
\leq c\, C_1 (CN)^{n} \sum_{\ell + m = n} 2^{2\ell }\sum_{i = 0}^\ell \frac{(2c)^i r_0^{2i} (n+1)^{i+1}}{i! }.
\end{multline*}
Therefore, there exist a constant $D$ depending on $r_0$ and $\Norm{P}{}$ such that 
$$
\forall\, n \geq N, \quad \Norm{K_n}{} \leq (DN)^{n}. 
$$

As $K_n$ is a polynomial of order at most $r_0 n$, we deduce from the previous estimate and Proposition \ref{P1} that, for $z \in B^s_{9M}$, 
$$
\Norm{X_{K_n}(z)}{\ell^1_s}\leq 2 (nr_0)^{s+1} (DN)^{n}(9M)^{ nr_0}\leq ( 2 c_s r_0^{s+1} DN(9M)^{ r_0})^n, 
$$
where the constant $c_s$ is defined in \eqref{Ecs}. Hence the series $\sum_{n\geq 0} t^{n}X_{K_n}(z)$ converges  for $t\leq ( 4 c_s r_0^{s+1} DN(9M)^{ r_0})^{-1}$.  
Furthermore, again for  $z\in B^s_{9M}$ and $t\leq ( 4 c_s r_0^{s+1} DN(9M)^{ r_0})^{-1}$, we  get using \eqref{Est1} and the previous bound
$$
\Norm{X_{Q_h(t)}(z) - X_P(z)}{\ell^1_s} \leq  \sum_{n\geq N} t^n \Norm{X_{K_n}(z)}{\ell^1_s} \leq (N+1)\, t^{N} (BN)^{N}
$$
for some constant $B$ depending   on $\Norm{P}{}$, $s$, $\delta$ and $r_0$. 

Let us  set
\begin{multline*}
c_0(M,r_0,\delta,s,\Norm{P}{}):= \\
 \min\{(4r_0^{s+1}|\Norm{P}{}M^{(r_0-1)})^{-1}, (4C c_s r_0^{s+1}c )^{-1}(9M)^{-r_0},( 4 c_s r_0^{s+1} D(9M)^{ r_0})^{-1}\}.
\end{multline*}
For $t\leq c_0(M,r_0,\delta,s,\Norm{P}{})N^{-1}$, inserting  the last  estimate  in 
\eqref{Est2} we get
$$
\Norm{u(t)}{\ell^1_s} \leq t^{N+1} (B N)^{N} +L_P \int_{0}^t \Norm{u(s)}{\ell^1_s} \dd s. 
$$
and this leads to
$$
\Norm{u(t)}{\ell^1_s} \leq t^{N+1} (\tilde{B}N)^{N}
$$
for some constant $\tilde{B}$ depending  on $r_0$, $\delta$, $s$  $\Norm{P}{}$ and $M$. 
This  implies \eqref{Estfond} defining $H_h=Z_h(h)/h$ for $h \leq c_0(M,r_0,\delta,s,\Norm{P}{})N^{-1}$.  

The second assertion of the theorem is just a calculus defining
$$
H_h^{(1)}=\frac{1}hZ_0+ Z_1.
$$
Using the previous bounds and the first inequality in Proposition \ref{P1}, we then calculate  that for $z \in B^s_M$
\begin{align*}
\Norm{H_h(z) - H_h^{(1)}(z)}{\ell^1_s} &\leq \sum_{j = 2}^{N} h^{j-1}\Norm{Z_j(z)}{\ell^1_s}\\[2ex]
&\leq  h c M^{2 r_0} \sum_{j = 2}^N h^{j-2}(C j M^{r_0})^{j-2}\\[2ex]
&\leq h c M^{2 r_0}\sum_{j = 2}^N \Big(\frac{j}{2N}\Big)^{j-2} \leq 2h M^{2 r_0}
\end{align*}
by definition  of $c_0(M,r_0,\delta,s,\Norm{P}{})$.  
\end{Proof}

\begin{Proofof}{Lemma \ref{Linvert}}
With the previous notations, we have 
$$
Q_{h}(t) = \sum_{n \geq 0} t^n K_n 
$$
where $K_n$ is given by \eqref{EKn} and the bounds obtained show that $Q_h(t)(z)$ and $X_{Q_h(t)}(z)$ are well defined on $B_{9M}$. Now let us consider the flow $\Phi_{Z_h(t)}^1$. As previously mentioned, it acts from $B^s_M$ to $B^s_{9M}$. Now we can write formally
\begin{equation}
\label{EPHIZ}
\begin{array}{rcl}
\Phi_{Z_h(t)}^1(z) &=& \displaystyle\sum_{k \geq 0} \frac{1}{k!} (\Lc_{Z_{h}(t)})^k\\[2ex]
&=& \displaystyle\sum_{n \geq 0} t^n \sum_{k \geq 0} \frac{1}{k!} \sum_{\ell_1 + \cdots + \ell_k = n| \ell_i \leq N} 
\Lc_{Z_{\ell_1}} \circ \cdots \circ \Lc_{Z_{\ell_l}}[I](z)\\[2ex]
&=& \displaystyle\sum_{n \geq 0} t^n \Psi_n(z).
\end{array}
\end{equation}
We are going to show that this series converges uniformly for  $z\in B^s_M$ and $t\leq T$. \\
Let $K$ be fixed polynomial of degree $k$, and $G(z) = (G_j(z))_{j \in \Zc}$ a function acting on $\ell^1_s$, and taking value in $\ell^1_s$, and such that the entries $G_j(z)$ are all polynomials of degree $\ell$. 
By definition, we have 
$$
(\Lc_K \circ G)_j =\displaystyle\sum_{i \in \Zc} (X_K)_i \frac{\partial G_j}{\partial z_i} = \displaystyle\{K,G_j\}. 
$$
Now using the relation \eqref{Ebrack} of Proposition \ref{P1}, we have that $(\Lc_K \circ G)_j$ is a polynomial of degree $k + \ell - 2$ and of norm
smaller that $2k\ell \Norm{K}{}\Norm{G_{j}}{}$. Now if $K = Z_0$, this bounds can be refined in $\Norm{(\Lc_{Z_0} \circ G)_j}{} \leq (2\pi - \delta)\Norm{G_j}{}$ using \eqref{Enonres3}. \\
For a given $z \in B^s_M$, the $j$-th component of $\Lc_{Z_{\ell_1}} \circ \cdots \circ \Lc_{Z_{\ell_l}}(I)(z)$ is a polynomial of order $n(r_0-2) + 2$ and involve terms of momentum $\M(\jb) = - \epsilon a$ if $j = (a,\epsilon)$ (see \eqref{Emomentj}). Hence summing in $j$, and separating as before the indices $m$ for which $\ell_m = 0$ in the sum, we  obtain for a given $n$ and $z \in B^s_M$, 
\begin{multline*}
\Norm{\Psi_n(z)}{\ell^1_s} \leq 2 (nr_0)^{s+1} M^{r_0n} \times \\
\sum_{i = 1}^n \sum_{k \geq i }\sum_{\ell_1 + \cdots + \ell_i = n| \ell_n \leq N}  \frac{(2\pi - \delta)^{k-i}}{(k-i)! i! }(2 n r_0^2)^i 
 \ell_1\Norm{Z_{\ell_1}}{} \cdots \ell_i\Norm{Z_{\ell_i}}. 
\end{multline*}
and we conclude  as before that this series is convergent for $t \leq T$ given in \eqref{EestT} with $K$ depending  on $r_0$, $M$, $\Norm{P}{}$, $s$ and $\delta$. 

Now writing down the same argument for the series (in $k \geq 0$ and $t^n$, $n \geq 0$) defining $\frac{\dd}{\dd t} \Phi_{Z_h(t)}^1$ and $X_{Q_h(t)}\circ \Phi_{Z_h(t)}^1$, we see that this series are again convergent, which justify the relation \eqref{Erelat}. 
\end{Proofof}

\section{Applications}

In this section we analyze the consequences of the  analytic estimates obtained in the previous section. We first show that for the ``new schemes'' in table 1, we obtain exponential estimates. We then show that for general splitting schemes, we obtain results under an additional CFL condition depending on $r$. 

\subsection{Exponential estimates}

We consider  the following splitting scheme
\begin{equation}
\label{Elie1}
\Phi_P^h \circ \Phi_{A_0}^1
\end{equation}
where the operator $A_0$ is associated with a function $\alpha_h(x)$ (see \eqref{Ealphah}) satisfying 
\begin{equation}
\label{Econda}
\forall\, x \in \R\, \quad |\alpha_h(x)| \leq \gamma h^\beta
\end{equation}
for some $\beta \in (0,1)$ and some $\gamma > 0$. Examples of such methods preserving the global order $1$ approximation property of the scheme \eqref{Elie1} for smooth functions are given in table 1. 
For such scheme, we obtain an exponentially closed modified energy:
\begin{theorem}
\label{Texp}
Let $r_0 \geq 3$  $s \geq 0$ and $M \geq 1$ be  fixed. Assume that $P \in \Pc_{r_0}$ and that $\alpha_h$ satisfies the condition \eqref{Econda} for some constants $\gamma$ and $\beta$. Then there exists a constant $h_0$ depending   on $r_0$, $\Norm{P}{}$, $s$, $M$ and $\gamma$ such  that for all $h\leq h_0$, there exists a real polynomial Hamiltonian $H_h $ such that for all $z \in B^s_M$, we have 
\begin{equation}
\label{Eexp}
\Norm{\Phi_P^h \circ \Phi_{A_0}^1(z) - \Phi_{H_h}^h(z) }{\ell^1_s} \leq h  \exp( - (h_0/h)^{\beta}). 
\end{equation}
\end{theorem}
\begin{Proof}
The hypothesis \eqref{Econda} implies that the eigenvalues $\lambda_a$ of the operator $A_0$ are bounded by $\gamma h^\beta$. Hence for a multi-index $\jb = (j_1,\ldots, j_r)$ we have 
$$
|\Lambda(\jb)| \leq r \gamma h^{\beta}
$$
and the condition $|\Lambda(\jb)| \leq \pi$ will be satisfied as long as $r \leq (\pi/\gamma) h^{-\beta}$. 
Taking $r_1$ such that $r_1 \leq (\pi/\gamma) h^{-\beta} < r_1 +1$ and defining $N = (r_1 - 1)/(r_0 - 1)$, we get 
$ b_1 h^{-\beta} \leq N\leq  b_2 h^{-\beta}$ for some positive constants $b_1$ and $b_2$ depending on $\gamma$ and $r_0$. 
Now the estimate \eqref{Estfond} in Theorem \ref{TB} for this $N$ yields
$$
\Norm{\Phi_P^h \circ \Phi_{A_0}^1(z) - \Phi_{H_h}^h(z) }{\ell^1_s} \leq h^{N+1} (CN)^N \leq h (Cb_2 h^{1- \beta})^{N},
$$
as long as $hN \leq c_0$. 
Thus  defining $$h_0 = \min\{(c_0b_1^{-1})^{1/(1 - \beta)},(eCb_2)^{- 1/(1 - \beta)}, b_1^{1/\beta}\},$$
we have for $0<h\leq h_0$ 
$$
\Norm{\Phi_P^h \circ \Phi_{A_0}^1(z) - \Phi_{H_h}^h(z) }{\ell^1_s}\leq h e^{-N} \leq  h e^{- b_1 h^{-\beta}}\leq h\exp( - (h_0/h)^{\beta}).
$$  
\end{Proof}

The dynamical consequences for the associated numerical scheme are given in the following corollaries.   We first assume that the numerical solution remains a priori in $\ell^1$ {\em only} over arbitrary long time.
\begin{corollary}
\label{ECC}
Under the hypothesis of the previous Theorem, let $z^0 = (\xi^0,\bar\xi^0) \in \ell^1$ and the sequence $z^n$ defined by 
\begin{equation}
\label{Ereczn}
z^{n+1} = \Phi_P^h \circ \Phi_{A_0}^1(z^n), \quad n \geq 0. 
\end{equation}
Assume that for all $n$, the numerical solution $z^n$ remains in a ball $B_M$ of $\ell^1$ for a given $M>0$. Then there exist constants $h_0$ and $c$  such that for all $h \leq h_0$, there exists a polynomial Hamiltonian $H_h$ such that 
$$
H_h(z^n) = H_h(z^0)  + \mathcal{O}(\exp(-ch^{-\beta}))
$$
for $nh \leq \exp(ch^{-\beta})$. Moreover, with the Hamiltonian $H_h^{(1)}$ defined in \eqref{EH1} we have 
\begin{equation}
\label{EpresH1}
H_h^{(1)}(z^n) = H_h^{(1)}(z^0)  + \mathcal{O}(h)
\end{equation}
over exponentially long time $nh \leq \exp(ch^{-\beta})$. 
\end{corollary}
This means that the modified energy remains exponentially closed  from its initial value during exponential times. Or more practically (since $H_h^{(1)}$ is explicit) the first modified energy is almost conserved during exponential times. 

\begin{Proof}
As all the Hamiltonian function considered are real, we have for all $n$, $z^n = (\xi^n,\bar\xi^n)$, i.e. $z^n$ is real. Hence for all $n$, $H_h(z^n) \in \R$.  Note moreover that we can always assume that $M \geq 1$.   \\
We use the notations of the previous theorems and we notice that  $H_h(z)$ is a conserved quantity by the flow generated by  $H_h$. Therefore we have
$$
H_h(z^{n+1}) - H_h(z^n) =  
H_h(\Phi_P^h \circ \Phi_{A_0}^1(z^n)) - H_h(\Phi_{H_h}^h(z^n)) 
$$
and hence
$$
|H_h(z^{n+1}) - H_h(z^n)| \leq \Big(\sup_{z \in B_{2M}}\Norm{\nabla H_h(z)}{\ell^\infty}\Big) \Norm{\Phi_P^h \circ \Phi_{A_0}^1(z^n) - \Phi_{H_h}^h(z^n)}{\ell^1}. 
$$
Now using \eqref{Egrad} and \eqref{Eexp}
we obtain for all $n$
$$
|H_h(z^{n+1}) - H_h(z^n) |\leq 4\pi M h  \exp( - (h_0/h)^{\beta})
$$
and hence
$$
|H_h(z^{n}) - H_h(z^0) |\leq (nh)  \exp( - 2c h^{-\beta})
$$
for some constant $c$, provided $h_0$ is small enough. 
This implies the result. The second estimate is then a clear consequence of \eqref{EestH1}.
\end{Proof}

The preservation of the Hamiltonian $H_h^{(1)}$ over long time induces that for $z^n = (\xi^n,\eta^n)$, 
$$
\sum_{a \in \Nc} \frac{1}{h}\alpha_h(h\omega_a) \xi^n_a\eta^n_a
$$
is bounded over long time, provided $z^0$ is smooth. For the functions $\alpha_h$ given in Table 1, this yields $H^{m/2}$ bounds for low modes, while the $L^2$ norm of high modes remain small (see \cite[Corollary 2.4]{DF09} for similar results in the linear case). We detail here the result in the specific situation where $\alpha_h(x)$ is given by \eqref{Egood}.
\begin{corollary}
\label{ECCC}
Let $r_0 \geq 3$,  $P \in \Pc_{r_0}$ and  $
\alpha_h(x) = \sqrt{h}\arctan(x/\sqrt{h}). $
We assume that  there exists a constant $b \geq 1$ such that 
\begin{equation}
\label{Edown}
\forall\, a \in \Nc,\quad  \frac{1}{b} |a|^m \leq \omega_a \leq b |a|^m. 
\end{equation}
Let $z^n$ the sequence defined by \eqref{Ereczn}. We assume  that $z^0 = (\xi^0,\bar\xi^0) \in H^{m/2}$ and that the sequence $z^n$ remains bounded in $\ell^1$. Then there exists constants $C$, $\alpha$ and $\beta$ such that 
$$
\sum_{|j| < \alpha h^{-1/2m}} |j|^{m} |z_j^n|^2 + \frac{1}{\sqrt{h}} \sum_{|j| \geq \alpha h^{-1/2m}}  |z_j^n|^2 \leq C
$$
over exponentially long time $nh \leq \exp(ch^{-\beta})$. 
\end{corollary}

\begin{Proof}
The almost conservation  of the Hamiltonian $H_h^{(1)}$ (cf. \eqref{EpresH1} ) shows that for all $n$ such that $nh \leq \exp(ch^{-\beta})$, we have
\begin{multline*}
\sum_{a \in \Nc} \frac{1}{\sqrt{h}}\arctan(\sqrt{h}\omega_a) |\xi^n_a|^2\\[2ex]
\leq |P_h^{(1)}(z^0) - P_h^{(1)}(z^n)| + \sum_{a \in \Nc} \frac{1}{\sqrt{h}}\arctan(\sqrt{h}\omega_a) |\xi^0_a|^2 + Ch
\end{multline*}
for some constant $C$, where $P_h^{(1)}(z)$ is the non linear Hamiltonian of equation \eqref{EH1}. 
As $z^n$ remains bounded in $\ell^1$, we have that $|P_h^{(1)}(z^0) - P_h^{(1)}(z^n)|$ is uniformly bounded. Now we have for all $a$ 
$$
\frac{1}{\sqrt{h}}\arctan(\sqrt{h}\omega_a)\leq \omega_a. 
$$
and we deduce using  \eqref{Edown}
$$
\sum_{a \in \Nc} \frac{1}{\sqrt{h}}\arctan(\sqrt{h}\omega_a) |\xi^0_a|^2  \leq C \Norm{z}{H^{m/2}}.
$$
Hence we have 
$$
\sum_{a \in \Nc} \frac{1}{\sqrt{h}}\arctan(\sqrt{h}\omega_a) |\xi^n_a|^2 \leq C
$$ 
for some constant $C$ depending on $M$ and $\Norm{z}{H^{m/2}}$. We conclude by using 
$$
x > 1 \Longrightarrow\arctan x > \arctan(1)\quad \mbox{and}\quad x \leq 1 \Longrightarrow \arctan x > \frac{x}{2}. 
$$
and the bounds \eqref{Edown} on $\omega_a$. 
\end{Proof}

Now we conclude this section by considering the case where  the numerical solution $z^n$ remains a priori bounded in some $\ell^1_s$, $s > 0$. In this case, we can show that the initial energy $H(z^n)$ remains bounded over exponentially long time, provided that $s$ is sufficiently large (typically of order 3 in concrete applications). We first begin with the following  Lemma: 

\begin{lemma}
\label{LCC}
Assume that $h \leq 1$, and that $\alpha_h(x)$ satisfies: 
\begin{equation}
\label{Econdalph2}
\forall\, x > 0, \quad  |\alpha_h(x) -x | \leq C h^{-\sigma} x^{\gamma},
\end{equation}
for some constants $\sigma \geq 0$, $\gamma \geq 2 + \sigma$ and $C> 0$. 
Assume moreover the non resonance condition \eqref{Enonres} is satisfied for some $\delta > 0$, and let $H_h^{(1)}$ be the Hamiltonian defined in \eqref{EH1}. Then for all $M > 0$ and $s \geq m\gamma/2$ (where $m$ is given in \eqref{Eboundomega}), and for all $z \in B^s_M$, 
$$
|H_h^{(1)}(z)  - H(z) | \leq C_1 h
$$
where $C_1$ depends on $M$ and $\delta$. 
\end{lemma}
\begin{Proof}
For all $a \in \Nc$, we have using \eqref{Econdalph2}   and $\gamma -\sigma - 1\geq1$ 
$$
\left| \frac{1}{h}\alpha_h(h\omega_a) - \omega_a \right| \leq C h^{\gamma -\sigma - 1} \omega_a^{\gamma},
 \leq Ch \omega_a^\gamma .
$$ 
Hence as $z = (\xi,\eta)$ is real, 
$$
\left|\sum_{a \in \Nc} \frac{1}{h}\alpha_h(h\omega_a) \xi_a\eta_a - H_0(z) \right| \leq C h \sum_{a \in \Nc} \omega_a^{\gamma} \xi_a\eta_a\leq C h \Norm{z}{\ell^1_{m\gamma/2}}^{2} \leq C M^2 h. 
$$
Now, for all $|x| \leq 2\pi - \delta$, we have 
$$
\left| 1 - \frac{ix}{e^{ix} - 1}\right| \leq C_\delta |e^{ix} - 1 - ix| \leq C_\delta |x|^2,
$$
for some constant $C_\delta$ depending on $\delta$. Hence we have for all $z \in \ell^1_s$, 
$$
\begin{array}{rcl}
\displaystyle\left|\sum_{\ell = 1}^{r_0} \sum_{\jb \in \Ic_\ell} \frac{i\Lambda(\jb)}{\exp(i \Lambda(\jb)) - 1} a_{\jb} z_{\jb} - P(z)\right|
&\leq& \displaystyle C_\delta \sum_{\ell = 1}^{r_0} \sum_{\jb \in \Ic_\ell}  |\Lambda(\jb)|^2 |z_{\jb}|\\[2ex]
&\leq& \displaystyle C_\delta (2\pi - \delta) \sum_{\ell = 1}^{r_0} \sum_{\jb \in \Ic_\ell}  |\Lambda(\jb)| |z_{\jb}|.
\end{array}
$$
Now  using the fact that for a given multiindex $\jb = (j_1,\ldots,j_{r_0})$ with for all $i = 1,\ldots,r_0$, $j_i = (a_i,\delta_i)$ , we have 
$$
|\Lambda(\jb)| \leq r_0 h \max_{i = 1,\ldots,r_0} |\omega_{a_i}| \leq r_0 h \prod_{i = 1}^{r_0} |j_i|^{m},
$$
we obtain 
$$
\displaystyle\left|\sum_{\ell = 1}^{r_0} \sum_{\jb \in \Ic_\ell} \frac{i\Lambda(\jb)}{\exp(i \Lambda(\jb)) - 1} a_{\jb} z_{\jb} - P(z)\right|
\leq C_\delta(2\pi - \delta ) r_0 h \Norm{z}{\ell^1_m}^{r_0}. 
$$
As $m\gamma/2 \geq m$, we have $\Norm{z}{\ell^1_m}^{r_0} \leq M^{r_0}$ , and this yields the result. 
\end{Proof}

As an example of application, we give the following result: 
\begin{corollary}
\label{CCC1}
Under the hypothesis of Corollary \ref{ECC}, 
assume that $\alpha_h(x) = \sqrt{h} \arctan(x/\sqrt{h})$, and that the sequence $z^n$ defined by \eqref{Ereczn} remains in a ball $B^s_M$ of $\ell^1_s$ for a given $M>0$ and $s \geq 3m/2$ where $m$ is given in \eqref{Eboundomega}. Then there exist constants $h_0$ and $c$  such that for all $h \leq h_0$, and all $n$ such that $nh \leq \exp(c/\sqrt{h})$, we have  
\begin{equation}
\label{EpresH}
H(z^n) = H(z^0)  + \mathcal{O}(h)
\end{equation}
where $H = H_0 + P$ is the original Hamiltonian.  
\end{corollary}

\begin{Proof}
In view of Corollary \ref{ECC}, we only have to prove that for real $z \in B^s_M$, we have 
$$
|H_h^{(1)}(z)  - H(z) | \leq C h
$$
which is given by the previous Lemma with $\gamma = 3$,  $\sigma = 1$ and $\delta = \pi$. 
\end{Proof}

\subsection{Results under CFL condition}

We now consider cases where $\alpha_h$ does not depend on $h$ and thus does not satisfy \eqref{Econda}. 
We focus on schemes such that 
\begin{equation}
\label{EalphCFL}
\alpha_h(x) = \beta(x) \mathds{1}_{x < c}(x)
\end{equation}
i.e. schemes associated with a filter function $\beta$ and a CFL condition with CFL number $c$. We mainly have in mind the two applications (see Table 1)
$$
\beta(x) = x \quad \mbox{and} \quad \beta(x) = 2 \arctan(x/2)
$$
corresponding to standard and implicit-explicit splitting schemes. 
Now for such scheme and for all multi-index $\jb = (j_1,\ldots, j_r)$ we have at most
$$
|\Lambda(\jb)| \leq r \beta(c). 
$$
Hence if we define
\begin{equation}
\label{Ecr}
c_r = \beta^{-1}(\frac{2\pi}{r}),
\end{equation}
the condition \eqref{Enonres} will be satisfied for some $\delta$ if $c < c_r$. The following result is then a easy consequence of Theorem \ref{TB}: 
\begin{theorem}
\label{TCFL}
Let $r_0 \geq 3$, $r \geq r_0$, $s \geq 0$ and $M >0$ be fixed. Assume that $P \in \Pc_{r_0}$ and that $\alpha_h$ is of the form \eqref{EalphCFL} for some constants $c < c_r$ where  $c_r$ is given by the equation \eqref{Ecr}. Then there exist constants $h_0$  and $C_r$ such that  for all $h\leq h_0$, there exists a real polynomial Hamiltonian $H_h$ such that for all $z \in B^s_M$, we have 
\begin{equation}
\label{Eecfl}
\Norm{\Phi_P^h \circ \Phi_{A_0}^1(z) - \Phi_{H_h}^h(z) }{\ell^1_s} \leq C_r h^{N(r)}
\end{equation}
where $N(r) = (r-2)/(r_0 - 2) + 1$.\\
 If moreover $P \in \Sc\Pc_{r_0}$ and for all $a \in \Nc$, $\omega_a \geq 0$, then the same results holds true under the weaker condition
\begin{equation}
\label{Ebsym}
c <  c_{r/2} = \beta^{-1}(\frac{4\pi}{r}).
\end{equation}
\end{theorem}
\begin{Proof}
The first part of the Theorem is a consequence of the previous estimates. 
The last assertion comes from the fact that all the $\lambda_a$ are positive. Hence for a given $r$ and a given monomial of $\Sc\Pc_r$ associated with a symmetric multi-index $\jb = (j_1,\ldots,j_{r/2}, k_1,\ldots,k_{r/2})$ with $j_i = (a_i,+1)$ and $k_i = (b_i,-1)$, $a_i$ and $b_i \in \Nc$, we have 
$$
\lambda(\jb) = \lambda_{a_1} + \ldots + \lambda_{a_{r/2}} - \lambda_{b_1} - \ldots - \lambda_{b_{r/2}} \leq \frac{r\beta(c)}{2}
$$
and this yields the bound \eqref{Ebsym}. 
\end{Proof}

This result implies the preservation of the Hamiltonian over long times of order $h^{-N(r) + 1}$, as in Corollary \ref{ECC}. Similarly, Corollary \ref{ECCC} extends to the case of the mid-split scheme (i.e. the case $\beta(x) = 2 \arctan(x/2)$) leading to a $H^{m/2}$ control of the frequencies smaller than $c_rh^{-1/m}$, over long times of order $h^{-N(r) + 1}$,   namely
$$
\sum_{|j| < c_rh^{-1/m}} |j|^{m} |z_j^n|^2  \leq C
$$
for $n\leq h^{-N(r)}$.  As the proof is completely similar, we do not give the details here. Note however that the high frequency cut-off does not give a control of the $L^2$ norm for high modes. 

In Table 2, we give the expression of the CFL constant $c_r$ in \eqref{Ecr} required in order to obtain a given precision of order $h^{N(r)}$ in the estimate \eqref{Eecfl} in the two cases where $\beta(x) = x$ and $\beta(x) = 2\arctan(x/2)$. The third column is concerned with  the nonlinear cubic Schr\"odinger equation with a general cubic nonlinearity (i.e.\ a quartic Hamiltonian, $r_0 = 4$) while, in the fourth one, we consider  the case of the nonlinear cubic Schr\"odinger equation with a nonlinearity belonging to $\Sc\Pc_{4}$. For instance, we see that in this latter case, the preservation of the energy is  ensured over long times of order $h^{-6}$ with a CFL or oder $1$ for the mid-split integrator. 
\begin{table}[ht]
\begin{center} 
\begin{tabular}{|@{$\quad$}c@{$\quad$}|@{$\quad$}c@{$\quad$}|@{$\quad$}c@{$\quad$}|@{$\quad$}c@{$\quad$}|}
\hline
$h^{N(r)}$ & $\beta(x) = x$ &  $\beta(x) = 2\arctan(x/2)$ & $\begin{array}{c}\beta(x) = 2\arctan(x/2) \\ \mbox{for cubic NLS}\end{array}$ \\
\hline 
$h^2$ & $1.57$  &  $2.00$ & $\infty$ \\
\hline 
$h^3$ & $1.05$ & $1.15$  & $3.46$ \\ 
\hline 
$h^4$ & $0.79$ &  $0.83$ & $2.00$ \\    
\hline 
$h^5$ & $0.63$ &  $0.65$ & $1.45$ \\
\hline 
$h^6$ & $0.52$ &  $0.54$ & $1.15$ \\
\hline 
$h^7$ & $0.45$ &  $0.46$ & $0.96$ \\
\hline
$h^8$ & $0.39$ &  $0.40$ & $0.83$ \\
\hline
$h^9$ & $0.35$ &  $0.35$ & $0.73$ \\
\hline
$h^{10}$ & $0.31$ &  $0.32$ & $0.65$ \\
\hline 
\end{tabular}
\caption{CFL conditions for quartic nonlinearities}
\label{Table2}
\end{center}
\end{table}

\begin{commentary}
Note that Corollary \ref{CCC1} can be easily extended to the situation where $\beta(x) = 2 \arctan(x/2)$ with a CFL condition $c_r$. Indeed, in this case, for a real $z = (\xi,\eta)$ the difference between the energy $H_h^{(1)}(z)$ with $\alpha_h(x) = \beta(x) \mathds{1}_{x < c_r}(x)$ and the same energy but associated with $\alpha_h(x) = \beta(x)$ is 
$$
\sum_{a \in \Nc | h\omega_a \geq c_r} 2 \arctan(h\omega_a/2) \xi_a \eta_a \leq \pi \sum_{a \in \Nc | h|a|^m  \geq c_r} \xi_a \eta_a
$$ 
and this last expression is bounded by $C h \Norm{z}{\ell^1_s}^2$ for $s = m/2$. 
We then conclude using Lemma \ref{LCC} with $\alpha_h(x) = \beta(x)$, $\sigma = 0$ and $\gamma = 3$ to prove statements similar to Lemma \ref{CCC1} in this situation. 
\end{commentary}
%




%

\section{Equations on $\R^d$}

In this section we generalize our method to the ``continuous'' case, i.e.\ the case where the linear part of the Hamiltonian has a continuous spectrum. As an example we consider the nonlinear Schr\"odinger equation (NLS) on $\R^d$:
\begin{equation}\label{nls}
iu_t=-\Delta u +\partial_2 g(u,\bar u)
\end{equation}
where we assume that $g:\C^2\to\C$ is smooth and has a zero of order at least 3 at the origin and satisfies $g(z,\bar z)\in \R$.

The natural phase spaces are the Sobolev spaces, $(u,\bar u)\in H^s(\R^d)\times H^s(\R^d)$ which are endow with the symplectic structure associated with the following Poisson bracket
$$
\{F,G\}(u,\bar u)=i\int_{\R^d}\left(\frac{\partial F}{\partial u}\frac{\partial G}{\partial \bar u}-\frac{\partial F}{\partial \bar u}\frac{\partial G}{\partial u}\right) \dd x.
$$
In this Hamiltonian setting the NLS equation reads
$$
iu_t = \frac{\partial H}{\partial \bar u}, \quad i\bar u_t=-\frac{\partial H}{\partial u}$$
with
$$
H(u,\bar u) =\int _{\R^d}\left(|\nabla u(x)|^2 +g(u(x),\bar u(x))\right) \dd x.$$
The Fourier transform maps the symplectic space $H^s(\R^d)\times H^s(\R^d)$ into the symplectic space 
$\Pcc_s:=L^2_s(\R^d)\times L^2_s(\R^d)$ where
$$
L^2_s(\R^d):=\{\phi\in L^2_s(\R^d), \int_{\R^d}(1+|\xi|^{2s})|\phi(\xi)|^2\dd\xi < +\infty\}
$$
and, in the Fourier variables, the Poisson bracket reads
$$
\{F,G\}(\phi,\psi)=i\int_{\R^d}\left(\frac{\partial F}{\partial \phi}\frac{\partial G}{\partial \psi}-\frac{\partial F}{\partial \psi}\frac{\partial G}{\partial \phi}\right) \dd\xi.
$$
In this setting the NLS reads
$$
\left\{
\begin{array}{rcll}
i\phi_t &= \displaystyle\frac{\partial H}{\partial \psi}&=\xi^2+\displaystyle\frac{\partial P}{\partial \psi}, \\[2ex] i\psi_t&=-\displaystyle\frac{\partial H}{\partial \phi}&=-\xi^2-\displaystyle\frac{\partial P}{\partial \phi},
\end{array}
\right.
$$
with
$$
H(\phi,\psi) =\int _{\R^d}\xi^2\phi(\xi)\psi(\xi)dx +P(\phi,\psi)$$
and
$$
P(\phi,\psi)=\int_{\R^d}g\left(\int_{\R^d}\phi(\xi)e^{-ix\cdot \xi}\dd\xi, \int_{\R^d}\psi(\eta)e^{ix\cdot \eta}\dd\eta\right)\dd x.$$
In particular for the cubic NLS, $g(u,\bar u)=|u|^4$ and
$$
P(\phi,\psi)=\int_{\xi_1+\xi_2-\eta_1-\eta_2=0}\phi(\xi_1)\phi(\xi_2)\psi(\eta_1)\psi(\eta_2)\dd\xi_1\dd\xi_2\dd\eta_1\dd\eta_2.$$
More generally, i.e. for general (analytic) nonlinearity $g$,
\begin{multline*}
P(\phi,\psi)=\\[2ex]
\sum_{3\leq n}\sum_{j+\ell=n}a_{j,\ell}\int_{\mathcal M(\xi,\eta)=0}\phi(\xi_1)\cdots\phi(\xi_j)\psi(\eta_1)\ldots \psi(\eta_\ell) \dd\xi_1\cdots \dd\xi_j\dd\eta_1\cdots \dd\eta_\ell
\end{multline*}
where 
$$
\mathcal M(\xi,\eta)=\xi_1+\cdots+\xi_j-\eta_1-\cdots-\eta_\ell$$ and
$$
a_{j,\ell}=\frac 1 {j!\ell !}\partial_1^j\partial_2^\ell g(0,0).
$$
To give a general framework similar to the discrete case 
we introduce the following notations: for a given $d$, we set 
$$
\Zc^c = \R^d \times \{\pm 1\}
$$
endowed with the canonical measure $\dd \xi \times \mathds{1}_{\{\pm 1\}}$.
For $\tau = (\xi,\delta)\in \Zc^c$,  we set $|\tau| = |\xi|$, and  we define the variable $z(\tau)=(\phi(\xi),\psi(\xi))$ by the formula
$$
 \left\{\begin{array}{rcll}
z( \xi,\delta) &=& \phi(\xi) & \mbox{if}\quad \delta = 1,\\[2ex]
z( \xi,\delta) &=& \psi(\xi) & \mbox{if} \quad \delta = -1. 
\end{array}\right.
$$
 For $s \geq 0$, 
we say that $z$ belongs to $L^1_s$  if the norm
$$
\Norm{z}{L_s^1} = \frac12 \int_{\R^d} (1 + |\xi|^s) (|\phi(\xi)| + |\psi(\xi)|) \dd \xi
$$
is finite.\\
Following Definition \ref{def:2.1}, we define
\begin{definition}\label{def:2.2}
We denote by $\Hc_s^c$ the space of Hamiltonians $P$ satisfying 
$$
P \in \mathcal{C}^{\infty}(L_s^1,\C), \quad \mbox{and}\quad 
X_P \in \mathcal{C}^{\infty}(L_s^1,L_s^1). 
$$
\end{definition}
For instance $$P(z)=\int_{\xi_1+\xi_2-\eta_1-\eta_2=0}\phi(\xi_1)\phi(\xi_2)\psi(\eta_1)\psi(\eta_2)\dd\xi_1\dd\xi_2\dd\eta_1\dd\eta_2$$ is in $\Hc_s^c$ with
$$X_P(z)(\xi,1)=2i\int_{\xi_1+\xi_2-\eta=\xi}\phi(\xi_1)\phi(\xi_2)\psi(\eta)\dd\xi_1\dd\xi_2\dd\eta$$
and
$$
 X_P(z)(\xi,-1)=-2i\int_{\xi_1-\eta_1-\eta_2=\xi}\phi(\xi_1)\psi(\eta_1)\psi(\eta_2)\dd\xi_1\dd\eta_1\dd\eta_2.
$$
As in the discrete case, we say that $z$ is real if $z = (\phi,\bar\phi)$ and that a Hamiltonian $P(z)$ is real if it is real for real $z$ (like the example just above). 

For the multi-variable $\taub = (\tau_1,\ldots,\tau_\ell) \in (\Zc^c)^{\ell}$ we use the notation
$$
z(\taub) = z(\tau_1)\cdots z(\tau_\ell). 
$$
We define the momentum $\Mc (\taub)$ of  the multi-variable $\taub$ by
$$
\Mc(\taub) = \xi_{1} \delta_{1} + \cdots + \xi_\ell \delta_\ell.
$$
We then define the subset of $(\Zc^c)^\ell$ 
$$
 \Ic_\ell^c =  \{  \taub = (\tau_1,\ldots,\tau_\ell) \in \Zc^{\ell}, \quad \mbox{with}\quad \Mc(\taub) = 0\}.
$$
We also use the notation $\taub = (\xib,\deltab)$ with $\xib = (\xi_1,\ldots,\xi_\ell)$ and $\deltab = (\delta_1,\ldots,\delta_\ell)$ which means that for all $i = 1,\ldots,\ell$ we have $\tau_{i} = (\xi_i,\delta_i)$. 
We denote by $\dd \Ic_\ell^c(\taub)$ the measure defined by 
$$
\int_{\Ic_\ell^c}f(\taub)\dd \Ic_\ell^c(\taub) = \sum_{\deltab \in \{ \pm 1\}^\ell }\int_{\Mc(\xib,\deltab) = 0} f(\xib,\deltab) \dd \xib
$$
where $\dd \xib = \dd \xi_1\cdots \dd \xi_{\ell}$.

\begin{definition}
We say that a polynomial Hamiltonian $P \in \Pc_{k}^c$ if it is real and if we can write 
$$
P(z) = \sum_{\ell = 2}^k \int_{\Ic_\ell} a(\taub) z(\taub) \dd \Ic_\ell^c(\taub)
$$
with 
$$
\sum_{\ell = 2}^{k}\Norm{a}{L^\infty(\Ic_\ell)} =: \Norm{P}{}< \infty. 
$$
\end{definition}
Similarly to the discrete case, we define the set of symmetric polynomials $\Sc\Pc_k$. 
The following result is the analog of Proposition \ref{P1}:

\begin{proposition}
Let $k \geq 2$  and $s \geq 0$, then we have $\Pc_k^c \subset \Hc_s^c$, and for $P \in \Pc_k^c$,  we have the estimates
$$
|P(z)| \leq \Norm{P}{}\big(\max_{n = 2,\ldots,k}  \Norm{z}{L^1_s}^n \big)
$$
and  for all $z \in L^1_s$,
$$
 \Norm{X_P(z) }{\ell^1_s} \leq 2 k (k-1)^{s} \Norm{P}{}\Norm{z}{L^1_s}\big(\max_{n = 1,\ldots,k-2}  \Norm{z}{L^1}^n \big). 
$$
Moreover, for $z$ and $y$ in $L^1_s$, we have 
$$
\Norm{X_P(z) - X_P(y) }{L^1_s} \leq 4 k (k-1)^s \Norm{P}{}{} \big(\max_{n = 1,\ldots,k-2} (\Norm{y}{L^1_s}^n, \Norm{z}{L^1_s}^n )\big)\Norm{z-y}{L^1_s}.
$$
Eventually, for $P\in \Pc_k^c$ and $Q \in \Pc_\ell^c$, then $\{P,Q\} \in \Pc_{k+ \ell - 2}^c$ and $$
\Norm{\{P,Q\}}{} \leq 2 k\ell\Norm{P}{}\Norm{Q}{}
$$
and if $P\in \Sc\Pc_k^c$ and $Q \in \Sc\Pc_\ell^c$, then $\{P,Q\} \in \Sc\Pc_{k+ \ell - 2}^c$.
\end{proposition}
\begin{Proof}
The first estimate is a clear consequence of the definition of $\Norm{P}{}$. 
Assume that $P$ is a homogeneous polynomial of degree $k$ with coefficients $a(\taub)$, $\taub \in \Ic_k$. Assume that $\taub \in \Ic_k$ is fixed. 
The derivative of a given monomial $z(\taub) = z(\tau_1)\cdots z(\tau_k)$ with respect to $z(\tau)$ vanishes except if $\tau \subset \taub$. Assume for instance that $\tau = \tau_k = (\zeta,\epsilon)$. Then the momentum zero condition implies that  $\Mc(\tau_1,\ldots,\tau_{k-1}) = -\epsilon \zeta$.  Hence  we can write with the previous notations
$$
(1 + |\tau|)^s
\left|\frac{\partial P}{\partial z(\tau)} \right|  \leq k  \Norm{P}{} 
\sum_{\deltab \in \{\pm 1\}^{k-1}} \int_{\xib \in (\R^d)^{k-1},\,  \Mc(\xib,\deltab) = -\epsilon\zeta}(1 + |\tau|)^s |z(\xib,\deltab)| \dd \xib.
$$
Now as in \eqref{raslebol} we have using the zero momentum condition
$$
(1 + |\tau|)^s \leq (1 + | \xi_1| + \ldots + |\xi_{k-1}|)^s \leq (i-1)^s \max_{n = 1,\ldots,k-1} (1 + |\xi_n|)^s
$$
As in \eqref{takecare}, we obtain the first estimate after  integrating in  $\tau =  (\zeta,\epsilon)$. 
The second estimate as well as the estimate on the Poisson bracket is proven similarly. 
\end{Proof}

We then consider operators $H_0$ with continuous frequencies $\omega(\xi)$. In the case of NLS  we have $\omega(\xi) = |\xi|^2$. 

Starting with this formalism, the formal analysis of section \ref{formal} is essentially the same with $A_0=\alpha_h(hH_0)$ for some real function $\alpha_h$  satisfying  $\alpha_h(x)\simeq x$ for small $x$. Thus
$$
A_0\phi(\xi)= \alpha_h(h\omega(\xi))\phi(\xi)\quad \mbox{and}\quad A_0\psi(\xi)= -\alpha_h(h\omega(\xi))\psi(\xi). 
$$
and we set $\lambda(\xi)=\alpha_h(h\omega(\xi))$. We consider the splitting method
$$
\Phi_P^h  \circ\Phi_{A_0}^1
$$
and we look for a real Hamiltonian function $Z(t,\phi,\psi)$ such that for all $t\leq h$ we have
$$
 \Phi_P^t  \circ\Phi_{A_0}^1=\Phi_{Z(t)}^1
$$
and $Z(0)= Z_0:=A_0$.

As in section \ref{formal}, plugging an Ansatz expansion $Z(t) = \sum_{\ell \geq 1} t^\ell Z_\ell$, we get the same recursive equations \eqref{Erec} to solve. 
The key lemma is the following: 

\begin{lemma}
Assume that 
$$
Q(z) = \int_{\Ic_k} a(\taub) z(\taub) \dd \Ic_k(\taub)
$$
is a homogeneous polynomial of order $k$. 
Then 
$$
\mathrm{ad}_{Z_0}(Q) = \int_{\Ic_k} i\Lambda(\taub) a(\taub) z(\taub) \dd \Ic_k(\taub)
$$
where for a multi-variable $\taub = (\tau_1,\ldots, \tau_r)$ such that for $i = 1,\ldots,r$, $\tau_i = (\xi_i,\delta_i) \in \R^d\times\{\pm 1\}$, we have 
$$
\Lambda(\taub) = \delta_1 \lambda(\xi_1) + \cdots + \delta_r \lambda(\xi_r). 
$$
\end{lemma}
Using this Lemma, we obtain that the condition \eqref{Enonres} takes the form
$$
\forall\, k = 3,\ldots,r,\quad
\forall\, \taub \in \Ic_k, \quad |\Lambda(\taub)| \leq 2\pi - \delta
$$
and the rest is extremely similar. For instance for the cubic NLS equation \eqref{nls} with $\alpha_h(x) = \sqrt{h}\arctan(x/\sqrt{h})$ and $g(u,\bar u)=|u|^4$   we get
\begin{theorem}
\label{t:nls}
Let $z^0 = (\phi_0,\bar\phi_0) \in L^1$ be the initial datum in Fourier space ($\phi_0=\hat u_0)$ and define the sequence $z^n=(\phi_n,\bar \phi_n)$  (the numerical solutions in Fourier space) by 
$$
z^{n+1} = \Phi_P^h \circ \Phi_{A_0}^1(z^n), \quad n \geq 0,
$$
where $A_0$ is associated with the filter function $\alpha_h(x) = \sqrt{h}\arctan(x/\sqrt{h})$ and $P$ is associated with the nonlinearity $\int |u|^4$.  
Assume that for all $n$, $z^n=(\phi_n,\bar \phi_n)$ remains bounded in $L^1$. Then there exist constants $h_0$ and $c$ such that for all $h \leq h_0$, there exists a polynomial Hamiltonian $H_h$ such that 
$$
H_h(\phi_n,\bar \phi_n) = H_h(\phi,\bar\phi)  + \mathcal{O}(\exp(-ch^{-1/2}))
$$
for $nh \leq \exp(ch^{-1/2})$.
Moreover, with the Hamiltonian $H_h^{(1)}$ given by
\begin{multline*}
H^{(1)}_h(\phi,\psi) = \int_{\R^d} \frac{1}{\sqrt h}\arctan(|\xi|^2/\sqrt{h})|\phi(\xi)|^2  \dd \xi + \\  \int_{\xi_1+\xi_2-\eta_1-\eta_2=0} \frac{i(|\xi_1|^2+|\xi_2|^2-|\eta_1|^2-|\eta_2|^2)}{e^{i (|\xi_1|^2+|\xi_2|^2-|\eta_1|^2-|\eta_2|^2)} - 1}\times \\   \phi(\xi_1)\phi(\xi_2)\psi(\eta_1)\psi(\eta_2)\dd\xi_1\dd\xi
_2\dd\eta_1\dd\eta_2
 \end{multline*}
 we have 
$$
H_h^{(1)}(\phi_n,\bar \phi_n) = H_h^{(1)}(\phi,\bar\phi)  + \mathcal{O}(h)
$$
over exponentially long time $n \leq \exp(ch^{-1/2})$. \\
Moreover, suppose that $u_0\in H^1(\R^3)$ then there exists a constant $C$ such that
$$
\int_{|\xi|\leq \alpha h^{-1/4}} |\xi|^2|\phi_n(\xi)|^2\dd\xi+ \frac 1 {\sqrt{h}} \int_{|\xi|\geq \alpha h^{-1/4}}|\phi_n(\xi)|^2\dd\xi \leq C 
$$
over exponentially long time $n \leq \exp(ch^{-1/2})$. \\
Finally, assume that for all $n$, $z^n$ remains in $B^s_M$ with $s \geq 3$, then we have 
$$
H(z^n) = H(z^0) + \mathcal{O}(h) 
$$
over exponentially long time $n \leq \exp(ch^{-1/2})$. 
\end{theorem}

\end{document}